\documentclass[12pt,leqno]{amsart}

\usepackage{amssymb,amsmath,amsthm,amscd}
\usepackage{mathrsfs}

\usepackage{enumerate}

\usepackage[all]{xy}
\CompileMatrices

\numberwithin{equation}{section}

\theoremstyle{plain}
\newtheorem{lemma}{Lemma}[section]
\newtheorem{prop}[lemma]{Proposition}
\newtheorem{theorem}[lemma]{Theorem}
\newtheorem{sublemma}[lemma]{Sublemma}
\newcommand{\Prop}{\begin{prop}}
\newcommand{\enprop}{\end{prop}}
\newcommand{\Lemma}{\begin{lemma}}
\newcommand{\enlemma}{\end{lemma}}
\newcommand{\Theorem}{\begin{theorem}}
\newcommand{\entheorem}{\end{theorem}}
\newtheorem{corollary}[lemma]{Corollary}
\newcommand{\Cor}{\begin{corollary}}
\newcommand{\encor}{\end{corollary}}
\newcommand{\Sub}{\begin{sublemma}}
\newcommand{\ensub}{\end{sublemma}}
\newtheorem{definition}[lemma]{Definition}

\newcommand{\Def}{\begin{definition}}
\newcommand{\edf}{\end{definition}}

\theoremstyle{definition}
\newtheorem{remark}[lemma]{Remark}

\newtheorem{conjecture}[lemma]{Conjecture}

\newcommand{\al}{\alpha}

\newcommand{\Aut}{\operatorname{Aut}}
\newcommand{\ba}{\begin{array}}
\newcommand{\bA}{{\mathbb A}}
\newcommand{\bbP}{\mathbb{P}}
\newcommand{\be}{\begin{enumerate}}

\newcommand{\bl}{\bigl}
\newcommand{\bnum}{\begin{enumerate}[{\rm (i)}]}
\newcommand{\bpf}{\begin{proof}}
\newcommand{\br}{\bigr}

\newcommand{\ch}{\on{ch}}
\newcommand{\cl}{\colon}
\newcommand{\cla}{\mathrm{cl}}

\newcommand{\codim}{{\operatorname{codim}}}

\newcommand{\C}{{\mathbb C}}

\newcommand{\Coh}{\on{Coh}}

\newcommand{\de}{\delta}

\newcommand{\dX}{{\mathop{X}\limits^\circ}}

\newcommand{\e}{{\mspace{1mu}\mathrm{e}\mspace{1mu}}}
\newcommand{\ea}{\end{array}}
\newcommand{\ee}{\end{enumerate}}

\newcommand{\eneq}{\end{eqnarray}}
\newcommand{\eneqn}{\end{eqnarray*}}
\newcommand{\eneqsub}{\end{eqnarray}\end{subequations}}
\newcommand{\enum}{\end{enumerate}}

\newcommand{\eq}{\begin{eqnarray}}
\newcommand{\eqn}{\begin{eqnarray*}}
\newcommand{\eqsub}{\begin{subequations}\begin{eqnarray}}
\newcommand{\ext}{{\mathscr{E}xt}}
\newcommand{\E}{\mathscr{E}}
\newcommand{\End}{\operatorname{End}}
\newcommand{\Ext}{\operatorname{Ext}}
\newcommand{\F}{\mathscr{F}}
\newcommand{\g}{{\mathfrak{g}}}
\newcommand{\gb}{\mathfrak{b}}
\newcommand{\gn}{\mathfrak{n}}
\newcommand{\gt}{\mbox{$\mspace{1mu}\mathfrak{t}\mspace{1.5mu}$}}

\newcommand{\Gro}{\mathcal{G}}
\newcommand{\h}{\mathscr{H}}

\newcommand{\hs}{\hspace}

\newcommand{\Hom}{\operatorname{Hom}}

\newcommand{\Ker}{\on{Ker}}
\newcommand{\la}{\lambda}
\newcommand{\lan}{\langle}
\newcommand{\La}{\Lambda}

\newcommand{\mm}{{\mathrm{E}\mspace{1mu}}}
\newcommand{\mono}{\rightarrowtail}
\newcommand{\monoto}{\rightarrowtail}

\newcommand{\ol}{\overline}
\newcommand{\on}{\operatorname}

\newcommand{\oplusl}{\mathop\oplus\limits}

\newcommand{\prolim}{\mathop{\varprojlim}\limits}
\newcommand{\pt}{\operatorname{pt}}

\newcommand{\Proof}{\begin{proof}}

\newcommand{\Q}{{\mathbb Q}}
\newcommand{\QED}{\end{proof}}
\newcommand{\ran}{\rangle}

\newcommand{\rf}{{\rm)}}

\newcommand{\ro}{{\rm(}}
\newcommand{\R}{{\rm R}}

\newcommand{\scup}{\mathop{\mbox{\small$\bigcup$}}}

\newcommand{\sD}{\mathscr{D}}

\newcommand{\set}[2]{\left\{#1\,;\,#2\,\right\}}
\newcommand{\seteq}{\mathbin{:=}}

\newcommand{\soplus}{\mathop{\mbox{\small$\bigoplus$}}\limits}
\newcommand{\sR}{\widetilde{\O_{X_w}}}

\newcommand{\ssum}{\mathop{\mbox{\small$\sum$}}}
\newcommand{\suml}{\sum\limits}
\newcommand{\supp}{\operatorname{supp}}

\newcommand{\Spec}{{\operatorname{Spec}}}
\newcommand{\Supp}{\operatorname{Supp}}
\newcommand{\tens}{\mathop\otimes\limits}

\newcommand{\tQ}{\widetilde{Q}}

\newcommand{\To}[1][\phantom{aaaa}]{\xrightarrow{\;#1\;}}

\newcommand{\Z}{{\mathbb{Z}\mspace{1mu}}}

\renewcommand{\ge}{\geqslant}
\renewcommand{\hom}{\operatorname{\it \mathscr{H}\kern-.25em om}}
\newcommand{\Rhom}{\mathrm{R}\kern-.2em\hom}

\renewcommand{\le}{\leqslant}

\renewcommand{\O}{\mathscr{O}}
\renewcommand{\phi}{\varphi}
\renewcommand{\P}{\mathbb{P}}

\newcommand{\isoto}[1][]{\mathop{\xrightarrow[#1]%
{\rule{0pt}{.9ex}%
{\raisebox{-.35ex}[0ex][-.6ex]{$\mspace{1mu}\sim\mspace{2mu}$}}}}}
\newcommand{\rbx}{\raisebox}
\newcommand{\sbullet}{\mspace{2mu}\mbox{\scriptsize$\bullet$}\mspace{2mu}}

\begin{document}

\title[Equivariant K-theory of affine flag manifolds]
{Equivariant K-theory of affine flag manifolds
and affine Grothendieck polynomials}
\author{Masaki KASHIWARA}
\address[Masaki KASHIWARA]{Research Institute for Mathematical Sciences,
Kyoto University, Kyoto 606--8502, Japan
}
\thanks{
MS is partially supported by NSF DMS-0401012.}
%
\author{Mark SHIMOZONO}
\address[Mark SHIMOZONO]{460 McBryde Hall,
         Department of Mathematics,
         Virginia Tech,
         Blacksburg, VA 24061-0123, USA
}
\thanks{}

\begin{abstract}
We study the equivariant K-group of
the affine flag manifold with respect to the Borel group action.
We prove that the structure sheaf of
the (infinite-dimensional) Schubert variety in the K-group
is represented by a unique polynomial, which we call
the affine Grothendieck polynomial.
\end{abstract}

\maketitle

\section{Introduction}
Let $G$ be a simply connected semisimple group, $B$ its Borel
subgroup, and $X=G/B$ the flag manifold.
Its $B$-orbits are of the form $BwB/B$ for some $w$ in the Weyl group $W$.
Its closure $X_w$ is called the Schubert variety.
It is well-known that the
equivariant K-group $K_B(X)$, which is the Grothendieck group of
the abelian category of coherent $B$-equivariant $\O_X$-modules, is a free
$K_B(\pt)$-module with $\{[\O_{X_w}]\}_{w\in W}$ as a basis.
Note that
$K_B(\pt)$ is isomorphic to the group ring $\Z[P]$ of the weight
lattice $P$ of a maximal torus of $B$. On the other hand,
$K_B(X)\simeq K_{B\times B}(G)$ gives another structure of
$K_B(\pt)$-module on $K_B(X)$ and we have a morphism
$\Z[P]\otimes\Z[P]\simeq K_B(\pt)\otimes K_B(\pt)\to K_B(X)$, which
factors through a homomorphism \eq \Z[P]\otimes_{\Z[P]^W}\Z[P]\to
K_B(X) \eneq called the equivariant Borel map. Here $\Z[P]^W$ is the ring of
invariants with respect to the action of the Weyl group $W$.
It is also well-known that $\Z[P]\otimes_{\Z[P]^W}\Z[P]\to  K_B(X)$
is an isomorphism. An element in $\Z[P]\otimes\Z[P]$ whose image is
$[\O_{X_w}]$, is known as a double Grothendieck polynomial when
$G=SL(n)$ \cite{LS}.

\medskip
The purpose of this paper is to generalize these facts to the affine case.
Contrary to the finite-dimensional case, there are two kinds of flag
manifolds; the inductive limit of Schubert varieties $\ol{BwB/B}$,
each of which is a finite-dimensional projective variety
(see \cite{Kumar} and the references there), and an
infinite-dimensional scheme whose Schubert varieties
$\ol{BwB_-/B_-}$ are finite-codimensional subschemes. Here, $B_-$ is
the opposite Borel subgroup.
In \cite{KK,Kumar}, Kostant-Kumar considered the first flag manifold
and studied its equivariant cohomology and K-theory.

In this paper we use the latter flag
manifold, which is studied in \cite{K}. We take the affine flag
manifold $X=G/B_-$ (see \S\,\ref{sec:flag}). It is an
infinite-dimensional (not quasi-compact) scheme over $\C$. Its
$B$-orbits are parameterized by the elements of the Weyl group $W$.
Each $B$-orbit $\dX_w$ is a locally closed subscheme with finite
codimension. As a scheme it is isomorphic to
$\bA^\infty=\Spec(\C[x_1,x_2,\ldots])$. The flag manifold $X$ is a
union of $B$-stable quasi-compact open subsets $\Omega$. We define
$K_B(X)$ as the projective limit of $K_B(\Omega)$. Then we have
$K_B(X)\cong\prod_{w\in W}K_B(\pt)[\O_{X_w}]$.

Similarly to the finite-dimensional case, we have a homomorphism
\eq&&\Z[P]\otimes_{\Z[P]^W}\Z[P]\to K_B(X). \label{eq:bor}\eneq In
the affine case this morphism is injective but is not surjective;
not all $[\O_{X_w}]$ are in the image of this morphism. However, as
we shall see in this paper, $[\O_{X_w}]$ is in the image after a
localization.

\smallskip

Let $\delta$ be the generator of null roots and let $R$ be the
subring of $\Q(\e^\delta)$ generated by $\e^{\pm\delta}$ and
$(\e^{n\delta}-1)^{-1}$ ($n\not=0$). Tensoring $R$ with \eqref{eq:bor},
we have the morphism
\eq&&R\tens_{\Z[\e^{\pm\delta}]} \Z[P]\tens_{\Z[P]^W}\Z[P]\to
R\tens_{\Z[\e^{\pm\delta}]}K_B(X). \label{eq:morgamma} \eneq Our
main result is the following theorem.

\medskip\noindent
{\bf Theorem \ref{th:main}}. {\it For all $w\in W$, $[\O_{X_w}]\in
K_B(X)$, considered as an element of
$R\otimes_{\Z[\e^{\pm\delta}]}K_B(X)$, is in the image of
\eqref{eq:morgamma}.}

\smallskip

Note that we have
$\Z[P]^W\simeq\Z[P^W]$ in the affine case.

\noindent Roughly speaking, we call the element of
$R\otimes_{\Z[\e^{\pm\delta}]}$$ \Z[P]\otimes_{\Z[P^W]}\Z[P]$
corresponding to $[\O_{X_w}]$ the {\em affine Grothendieck
polynomial} (see Proposition~\ref{prop:main}).

In order to prove Theorem~\ref{th:main},
we use the following vanishing theorem of the first group cohomology.

\smallskip\noindent
{\bf Theorem~\ref{th:van}}.
{\it
\bnum
\item
If $\vert I\vert>2$ then
$H^1(W,\,R\otimes_{\Z[\e^{\pm\delta}]}\Z[P])=0$.
\item
For any affine Lie algebra $\g$,
$H^1\bigl(W,R\otimes_{\Z[\e^{\pm\delta}]}
(\hs{-1ex}\smash{\soplus_{\la\in P,\,\vert\lan c,\la\ran\vert<\kappa^*}}\Z\e^\la)\bigr)=0$,
where $\kappa^*$ is the dual Coxeter number.
\ee
}

\noindent Here $I$ is the index set of simple roots and $c$ is the
canonical central element of $\g$.

\smallskip
The plan of this paper is as follows. In \S\,\ref{sec:flag} we
review the flag manifold of Kac-Moody Lie algebras. In
\S\,\ref{sec:dem} we study the Demazure operators. We also give a
simple proof of the fact that the Schubert varieties are normal and
Cohen-Macaulay. This proof seems to be new even in the
finite-dimensional case. In \S\,\ref{sec:aff} we study the affine
flag manifolds. After the preparation in \S\,\ref{sec:van}, we prove
Theorem~\ref{th:van} in \S\,\ref{sec:prvan}. As its application, we
give in \S\,\ref{sec:main} the proof of Theorem~\ref{th:main}, the
existence of affine Grothendieck polynomials. In
\S\,\ref{sec:charac} we give the character formula of the global
cohomology groups of $\O_{X_w}(\la)$ using the affine Grothendieck
polynomials. In \S\,\ref{sec:eqcoh} we explain an analogous result
for the equivariant cohomology groups of the affine flag manifolds.
In \S\,\ref{sec:ex} we shall give examples of the affine
Grothendieck polynomials.

\section{Flag manifolds}
\label{sec:flag} Let us recall in this section the definition and
properties of the flag manifold of a symmetrizable Kac-Moody Lie
algebra following \cite{K}.

Let $(a_{ij})_{i,j\in I}$ be a symmetrizable generalized Cartan
matrix, $\g$ an associated Kac-Moody Lie algebra, and
$\gt$ its Cartan subalgebra. Let $\g=\gn\oplus\gt\oplus\gn_-$ be the
triangular decomposition and
$\g=\soplus\nolimits_{\al\in\gt^*}\g_{\al}$ the root decomposition.
Let $\Delta\seteq\set{\al\in\gt^*}{\g_\al\not=0}\setminus\{0\}$ be
the set of roots and $\Delta^\pm$ the set of positive and negative
roots.

Let $\{\al_i\}_{i\in I}$ be the set of simple roots in $\gt^*$ and
$\{h_i\}_{i\in I}$ the set of simple co-roots in $\gt$. Hence we
have $\lan h_i,\al_j\ran=a_{ij}$. Let us take an integral weight lattice
$P\subset\gt^*$. We assume the following conditions:
\eq&&\left\{\parbox{70ex}{ \bnum
\item $\al_i\in P$ for all $i\in I$,
\item $\{\al_i\}_{i\in I}$ is linearly independent,
\item $h_i\in P^*$ for all $i\in I$, where $P^*$ is the dual lattice
$\Hom(P,\Z)\subset\gt$,
\item there exists $\Lambda_i\in P$ such that
$\lan h_j,\Lambda_i\ran=\delta_{ij}$.
\enum}\right.\label{cond:weight}
\eneq

Let $T$ be the algebraic torus with $P$ as its character lattice.
Let $W$ be the Weyl group. It is the subgroup of
$\Aut(\gt^*)$ generated by the simple reflections $s_i$ ($i\in I$)
:$s_i(\la)=\la-\lan h_i,\la\ran \al_i$.
Let $U_\pm$ be the group scheme with $\gn_\pm$
as its Lie algebra. Let $B_\pm=T\times U_\pm$ be the Borel subgroup,
whose Lie algebra is $\gb_\pm=\gt\oplus\gn_\pm$.
For any $i\in I$, let us denote by $P_i^\pm$ the parabolic group
whose Lie algebra is $\gb_\pm\oplus\g_{\mp\al_i}$.
Then $P_i^\pm/B_\pm$ is isomorphic to the projective line $\bbP^1$.

Let $P^+\,\seteq\,\set{\la\in P}{\text{$\lan h_i,\la\ran\ge0$ $(i\in I)$}}$
be the set of dominant integral weights.
For $\la\in P^+$, let $V(\la)$ (resp.\ $V(-\la)$) be
the irreducible $\g$-module with highest weight $\la$
(resp.\ lowest weight $-\la$).
Then $A(\g)\seteq\soplus_{\la\in P^+}V(\la)\otimes V(-\la)$
has an algebra structure and we denote $\Spec(A(\g))$ by $G_\infty$.
The scheme $G_\infty$ contains a canonical point $e$ and is endowed
with a left action of $P_i$ and a right action of $P_i^-$.
The union
of $P_{i_1}\cdots P_{i_m}eP_{j_1}^-\cdots P_{j_m}^-\subset G_\infty$
is an open subset of $G_\infty$ and we denote it by $G$. Then $P_i$
and $P_i^-$ act freely on $G$.
The flag manifold $X$ is defined as the quotient $G/B_-$. It is a
separated (not quasi-compact in general) scheme over $\C$.
It is covered by affine open subsets isomorphic to $\bA^\infty
\seteq\Spec(\C[x_1,x_2,\ldots])$ (or $\bA^n$),
and its structure sheaf $\O_X$ is coherent.
Let $x_0\in X$ be the image of $e\in G$. Then for $w\in W$, $wx_0\in X$
has a sense. The set $X$ has a Bruhat decomposition
$X=\bigsqcup_{w\in W}\dX_w$, where $\dX_w$ is the locally closed
subscheme $Bw x_0$ of $X$. Let $X_w$ be the closure of $\dX_w$
endowed with the reduced scheme structure.
It is called the {\em
Schubert variety}.
It has codimension $\ell(w)$, the length of $w$,
and its structure sheaf $\O_{X_w}$ is coherent.
As a set we have $X_w=\bigsqcup_{x\ge w}\dX_x$.

\begin{remark}
\bnum
\item
For any $w$, $\ol{B_-wB_-/B_-}$ is a finite-dimensional projective
subscheme of $X$ and its union $\cup_{w\in W}\ol{B_-wB_-/B_-}$ is an
ind-scheme. This is another flag manifold which we do not use
here.
\item
For a regular dominant integral weight $\la$, set $\widehat{V}(-\la)
=\prod_{\mu\in P}V(-\la)_\mu$ where $V(-\la)_\mu$ is the
weight space of $V(-\la$) of weight $\mu$. Then
$\P(\widehat{V}(-\la))=(\widehat{V}(-\la)\setminus\{0\})/\C^\times$
has a scheme structure, and $P_i$ acts on $\widehat{V}(-\la)$ and
$\P(\widehat{V}(-\la))$. Then $X$ is embedded in
$\P(\widehat{V}(-\la))$ by $x_0\mapsto \overline{u_{-\la}}$, where
$\overline{u_{-\la}}$ is the line containing the lowest weight
vector $u_{-\la}$. \ee
\end{remark}

Let $S$ be a finite subset of $W$ such that $x\in S$ as soon as
$x\le y$ for some $y\in S$. Then $\Omega_S\seteq\scup_{w\in S}\dX_w$
is a $B$-stable quasi-compact open subset which coincides with
$\scup_{w\in S}wB x_0$. Conversely, any $B$-stable quasi-compact
open subset is of this form. Let $\Coh_B(\O_{\Omega_S})$ be the
abelian category of coherent $B$-equivariant $\O_{\Omega_S}$-modules
and let $K_B(\Omega_S)$ be the Grothendieck group of
$\Coh_B(\O_{\Omega_S})$.
For $\F\in\Coh_B(\O_{\Omega_S})$, let us denote by
$[\F]$ the corresponding element of $K_B(\Omega_S)$.
For $w\in W$, the sheaf $\O_{X_w}$ is a coherent $B$-equivariant $\O_X$-module
and gives an element $[\O_{X_w}]$ of $K_B(\Omega_S)$.
The equivariant $K$-group $K_B(\Omega_S)$ is a module over the ring
$K_B(\pt)$. For $\la\in P$, we denote by $\e^\la$ the element of
$K_B(\pt)$ represented by the one-dimensional representation of $B$
given by $B\to T\To[{\e^\la}]\C^\times$. By $P\ni\la\mapsto\e^\la$,
$K_B(\pt)$ is isomorphic to the group ring $\Z[P]$. We also denote
by $\e^\la$ the element of $\Z[P]$ corresponding to $\la\in P$.

For $w\in S$, $w x_0\in\Omega_S$ is a $T$-fixed point. It defines a
$T$-equivariant inclusion $i_w\cl\pt\to \Omega_S$. Since any
coherent $\O_{\Omega_S}$-module $\F$ has locally a finite resolution
by locally free modules of finite rank (see Lemma~\ref{lem:81}),
the $k$-th left derived
functor $L_ki_w^*\F$ vanishes for $k\gg0$. Hence we can define the
$K_T(\pt)$-linear homomorphism
$$i_w^*\cl K_B(\Omega_S) \to K_T(\Omega_S)\to K_T(\pt)\simeq K_B(\pt)$$
by $[\F]\mapsto \sum_{k=0}^\infty(-1)^k[L_ki_w^*\F]$.
Note that, similarly to the finite-dimensional case, we have
\eq&&i_x^*([\O_{X_w}])=
\begin{cases}
\prod_{\al\in\Delta^+\cap w\Delta^-}(1-\e^\al)&\text{if $x=w$,}\\
0&\text{unless $x\ge w$.}
\end{cases}\label{eq:iwX}
\eneq

\begin{remark} The function $W\to \Z[P]$ given by the equivariant localization
$x\mapsto i_x^*([\O_{X_w}])$, coincides with
the function $\psi^w$ in \cite{KK}.
\end{remark}

\Lemma \label{lem:Omega} $K_B(\Omega_S)$ is a free $K_B(\pt)$-module
with basis $\{[\O_{X_w}]\}_{w\in S}$. \enlemma
\Proof Let us argue by induction on the cardinality of $S$. Let
$w$ be a maximal element of $S$. Set $S'=S\setminus\{w\}$. Then we
have $\Omega_S=\Omega_{S'}\sqcup\dX_w$. Hence we have an exact
sequence
$$K_B(\dX_w)\to K_B(\Omega_S)\to K_B(\Omega_{S'})\to 0.$$
By induction $K_B(\Omega_{S'})$ is a free
$K_B(\pt)$-module with a basis $\{[\O_{X_x}]\}_{x\in S'}$. Also
$K_B(\dX_w)$ is a free  $K_B(\pt)$-module generated by $[\O_{X_w}]$.
Hence $K_B(\Omega_S)$ is generated by $\{[\O_{X_x}]\}_{x\in S}$. By
\eqref{eq:iwX} the image of $\{[\O_{X_x}]\}_{x\in S}$  under the
map
$$K_B(\Omega_S)\To[{\prod_{x\in S}i_x^*}]K_B(\pt)^{\prod S}$$
is linearly independent over $K_B(\pt)$.
Here $K_B(\pt)^{\prod S}$ is the product of the copies of $K_B(\pt)$
parameterized by elements of $S$.
\QED

\begin{remark}
For $\ell\in\Z_{\ge0}$, let $\gn_\ell$ be the direct sum of $\g_\alpha$ where
$\al=\sum_im_i\al_i\in \Delta^+$ ranges over positive roots such that
$\sum_im_i>\ell$. Then $\gn_\ell$ is an ideal of $\gn$.
Let $U_\ell$ be the normal subgroup of $U$ with $\gn_\ell$ as its Lie algebra.
Then for any $S$ as above,
$U_\ell$ acts on $\Omega_S$ freely for $\ell\gg0$,
and the quotient space $U_\ell\backslash\Omega_S$
is a finite-dimensional scheme.
Hence $\Omega_S$ is the projective limit of
$\{U_\ell\backslash\Omega_S\}_{\ell}$ and
$K_B(\Omega_S)$ is the inductive limit of
$\{K_{B/U_\ell}(U_\ell\backslash\Omega_S)\}_\ell$.
\end{remark}

We set $K_B(X)=\prolim_SK_B(\Omega_S)$.
Hence we have
$K_B(X)=\prod\limits_{w\in W}K_B(\pt)[\O_{X_w}]$.

For $w\in W$, the homomorphism
$i_w^*\cl K_B(\Omega_S)\to K_T(\pt)$
induces a homomorphism
$$i_w^*\cl K_B(X)\to K_T(\pt).$$
By \eqref{eq:iwX} they induce a monomorphism
\eq
\xymatrix@C=10ex{{K_B(X)\ }\ar@{>->}[r]_(.45){\prod_{w\in W}i_w^*}
&K_T(\pt)^{\prod W}}. \label{eq:embK} \eneq

For $i\in I$, set $X_i=G/P_i^-$. Then there is a canonical projection
$p_i\cl X\to X_i$ which is a $\bbP^1$-bundle. We have
$p_i(\dX_w)=p_i(\dX_{ws_i})$ and
$p_i^{-1}p_i(\dX_w)=\dX_w\sqcup\dX_{ws_i}$ for any $w\in W$, and we
have the $B$-orbit decomposition $X_i=\mathop\sqcup\limits_{w\in
W,\, ws_i>w}p_i(\dX_w)$. Similarly to $K_B(X)$, we define $K_B(X_i)$
as $\prolim_SK_B(p_i(\Omega_S))$. It is isomorphic to
$\hs{-3ex}\prod\limits_{w\in W,\, ws_i>w} K_B(\pt)[\O_{p_i(X_w)}]$.

There exist homomorphisms
${p_i}_*\cl K_B(X)\to K_B(X_i)$ and
$p_i^*\cl K_B(X_i)\to K_B(X)$,
defined by
$[\F]\mapsto \sum_{k=0}^\infty(-1)^k[R^k{p_i}_*\F]=[{p_i}_*\F]-[R^1{p_i}_*\F]$
and $[\E]\mapsto [p_i^*\E]$, respectively.

Any element $\la\in P$ induces a group homomorphism $B_-\to
T\xrightarrow{\e^\la}\C^\times$. Let $\O_X(\la)$ be the invertible
$\O_X$-module on $X=G/B_-$ induced by this character of $B_-$. Then
we have a homomorphism of abelian groups $\Z[P]\to K_B(X)$ given by
$\e^\la\mapsto [\O_X(\la)]$.
Note that, for $\la\in P^+$,
we have $\Gamma(X;\O_X(\la))\simeq V(\la)$
and $H^k(X;\O_X(\la))=0$ for $k\not=0$ (see \cite{K1}).

Similarly to the finite-dimensional case (\cite{D,Kumar}),
we have a commutative diagram
\[
\xymatrix{
\Z[P]\ar[r]\ar[d]_{\sD_i}&K_B(X)\ar[d]^{p_i^*\circ{p_i}_*}\\
\Z[P]\ar[r]&K_B(X).}
\]
Here $\sD_i$ is given by \eq
&&\sD_i(\e^\la)=\dfrac{\e^\la-\e^{s_i\la-\al_i}}{1-\e^{-\al_i}},
\label{eq:demz} \eneq and is called the {\em Demazure operator}.

The $K_B(\pt)$-module structure on $K_B(X)$ induces a $K_B(\pt)$-linear
homomorphism $K_B(\pt)\otimes \Z[P]\to K_B(X)$. Let $\Z[P]^W$ be the
ring of $W$-invariants of $\Z[P]$. Then $K_B(\pt)\simeq\Z[P]$ are
$\Z[P]^W$-algebras. The morphism $K_B(\pt)\otimes \Z[P]\to K_B(X)$
decomposes as the composition of $K_B(\pt)\otimes \Z[P]\to
K_B(\pt)\otimes_{\Z[P]^W}\Z[P]$ and $$ \beta\cl
K_B(\pt)\otimes_{\Z[P]^W}\Z[P]\to K_B(X).$$
It is sometimes called the {\em equivariant Borel map}.

\begin{remark}
By \cite{KK}, the equivariant Borel map
$K_B(\pt)\otimes_{\Z[P]^W}\Z[P]\to K_B(X)$ is an isomorphism
if $G$ is a finite-dimensional simply connected semisimple group.
\end{remark}

The composition
\eq K_B(\pt)\otimes_{\Z[P]^W}\Z[P]\To[{\,\beta\,}]
K_B(X)\To[{\,i_w^*\,}]K_B(\pt)\simeq\Z[P] \label{eq:iw}
\eneq %
is given by $a\otimes b\mapsto a\cdot(wb)$.

\section{Demazure operators}
\label{sec:dem}

Let $p_i\cl X\to X_i$ be the $\bbP^1$-bundle as in \S\,\ref{sec:flag}.
In this section we shall show
\eqn
&&p_i^*{p_i}_*([\O_{X_w}])=
\begin{cases}
[\O_{X_{ws_i}}]&\text{if $ws_i<w$,}\\
[\O_{X_w}]&\text{if $ws_i>w$.}\\
\end{cases}
\eneqn %
Note that if $ws_i>w$, then $X_w=p_i^{-1}p_i(X_w)$ and
$p_i(X_w)=p_i(X_{ws_i})$. Moreover, $\dX_{ws_i}\to
p_i(\dX_{ws_i})=p_i(\dX_w)$ is an isomorphism.

\Lemma\label{lem:Rvan}
We have $R^1{p_i}_*\O_{X_w}=0$ for all $w\in W$.
\enlemma
\Proof
Assume that $ws_i<w$.
Since $X_{ws_i}=p_i^{-1}p_i(X_w)$, we have
$\O_{X_{ws_i}}=p_i^*\O_{p_i(X_w)}$, which implies that
$R^1{p_i}_*\O_{X_{ws_i}}=0$.
Applying the right exact functor $R^1{p_i}_*$ to
the exact sequence
$\O_{X_{ws_i}}\to\O_{X_{w}}\to0$,
we obtain
$R^1{p_i}_*\O_{X_{w}}=0$.
\QED

As shown in \cite{K}, for any point $p\in X_w$, there exist an open
neighborhood $\Omega$ of $p$ and a closed subset $S$ of $\bA^{n}$
for some $n$ such that we have a commutative diagram
$$\xymatrix{
{\rbx{-1ex}[1.5ex][2ex]{$X_w\cap\Omega$}}\ar@{-}^{\sim}[r]
\ar@{^{(}->}[d]&{\rbx{-1ex}[1.5ex][2ex]{$S\times\bA^\infty$}}\ar@{^{(}->}[d]\\
\Omega\ar@{-}^(.4){\sim}[r]&{\bA^n\times\bA^\infty}. }$$

Hence various properties of $X_w$ (such as normality) make sense.

\Prop\label{prop:normal}
For any $w\in W$, we have
\bnum
\item
$X_w$ is normal. \label{con:nor1}
\item
${p_i}_*\O_{X_w}\simeq\O_{p_i(X_w)}$.\label{con:nor2}
\ee
\enprop
\Proof
Let us show
\eqref{con:nor1} and \eqref{con:nor2} by induction on the length
$\ell(w)$.
Note that when $ws_i>w$, \eqref{con:nor2} follows from
$X_w=p_i^{-1}p_i(X_w)$.

When $w=e$, \eqref{con:nor1} and \eqref{con:nor2} are obvious.
Assume that $\ell(w)>0$.

\medskip
Let us first show \eqref{con:nor2}.
We may assume $ws_i<w$.
We have a monomorphism
$j\cl\O_{p_i(X_w)}={p_i}_*\O_{X_{ws_i}}\mono {p_i}_*\O_{X_w}$,
which is an isomorphism on $p_i(\dX_w)$.
By the induction hypothesis,
$X_{ws_i}$ as well as $p_i(X_w)=p_i(X_{ws_i})$
is normal.
Hence
$j$ is globally an isomorphism since
${p_i}_*\O_{X_w}$ is a coherent $\O_{X_i}$-module.

\medskip
Next let us show \eqref{con:nor1}. Let $\sR$ be the normalization of
$\O_{X_w}$. Then $\sR$ is a coherent $\O_X$-module. Let $S$ be the
support of $\sR/\O_{X_w}$. Then $S$ is a $B$-stable closed subset
contained in $X_w\setminus \dX_w$. We shall show that $S$ is an
empty set.

Otherwise let $x\in W$ be a minimal element such that $X_x\subset S$. Then
$x>w$. Let us take $i\in I$ such that $xs_i<x$. Assume first $ws_i>w$.
Then $\O_{X_w}=p_i^*\O_{p_i(X_w)}$. Hence $\sR$ is the inverse image
of the normalization of $\O_{p_i(X_w)}$. Hence $S=p_i^{-1}p_i(S)$.
This contradicts $X_{xs_i}\not\subset S$.
Hence we have
$ws_i<w$.

We have monomorphisms
$\O_{p_i(X_w)}\monoto {p_i}_*\O_{X_w}\monoto
{p_i}_*\sR$.
By the induction hypothesis,
$X_{ws_i}$ as well as $p_i(X_w)=p_i(X_{ws_i})$
is normal.
Since ${p_i}_*\sR$ is a coherent $\O_{p_i(X_w)}$-module and
$\O_{p_i(X_w)}\to{p_i}_*\O_{X_w}\to {p_i}_*\sR$ are isomorphisms
on $p_i(\dX_w)$,
the normality of $p_i(X_w)$ implies that
$\O_{p_i(X_w)}\monoto {p_i}_*\sR$ is an isomorphism on $X_i$.
Hence we have isomorphisms
$\O_{p_i(X_w)}\isoto {p_i}_*\O_{X_w}\isoto
{p_i}_*\sR$.

We have an exact sequence
${p_i}_*\O_{X_w}\to
{p_i}_*\sR\to {p_i}_*(\sR/\O_{X_w})\to R^1{p_i}_*\O_{X_w}$.
Since $R^1{p_i}_*\O_{X_w}=0$ by Lemma~\ref{lem:Rvan},
we obtain
${p_i}_*(\sR/\O_{X_w})=0$.
On the other hand, since $S\to X_i$ is an isomorphism on $p_i(\dX_x)$,
the support of
${p_i}_*(\sR/\O_{X_w})$
contains $p_i(\dX_x)$, which is a contradiction.
\QED

\Cor\label{cor:dem}
For any $w\in W$, we have
$$
p_i^*{p_i}_*([\O_{X_w}])=
\begin{cases}
[\O_{X_{ws_i}}]&\text{if $ws_i<w$,}\\
[\O_{X_w}]&\text{if $ws_i>w$.}
\end{cases}
$$
\encor

\Prop\label{prop:cm}
 The module $\O_{X_w}$ is Cohen-Macaulay for any $w\in W$.
\enprop \Proof Since $X_w$ has codimension $\ell(w)$,
$\ext^k_{\O_X}(\O_{X_w},\O_X)=0$ for $k<\ell(w)$. Hence it is enough
to show that $\ext^k_{\O_X}(\O_{X_w},\O_X)=0$ for $k>\ell(w)$. We
shall prove it by induction on $\ell(w)$. When $w=e$, it is
obvious. Assume that $\ell(w)>0$.

Set \eq \F=\tau^{>\ell(w)}\Rhom_{\O_X}(\O_{X_w},\O_X), \eneq where
$\tau^{>\ell(w)}$ is the truncation functor (see {\em e.g.}{}
\cite{KS}). Let us set $S=\Supp(\F)$. Then $S$ is a $B$-stable
closed subset of $X_w$. Let $x\in W$ be a minimal element of
$\set{x\in W}{X_x\subset S}$. Let us take $i\in I$ such that
$xs_i<x$. If $ws_i>w$, then we have
$$\F\simeq{p_i}^*\tau^{>\ell(w)}\Rhom_{\O_{X_i}}(\O_{p_i(X_w)},\O_{X_i}),$$
which implies that $p_i^{-1}p_i(S)=S$. Hence $X_{xs_i}\subset S$,
which is a contradiction. Hence we have $ws_i<w$.

Let $\Omega_{X/X_i}$ be the relative canonical sheaf, which is
isomorphic to $\O_X(-\alpha_i)$. Then the Grothendieck-Serre duality
theorem says that \eq
&&\R{p_i}_*\Rhom_{\O_X}(\O_{X_w},\Omega_{X/X_i}[1]) \simeq
\Rhom_{\O_{X_i}}(\R{p_i}_*\O_{X_w},\O_{X_i}). \label{eq:isorr}
\eneq Since $\R{p_i}_*\O_{X_w}\simeq\O_{p_i(X_w)}$
by Lemma~\ref{lem:Rvan} and Proposition~\ref{prop:normal} and since the
induction hypothesis implies that $p_i(X_w)=p_i(X_{ws_i})$ is
Cohen-Macaulay, we have
$$\Rhom_{\O_{X_i}}(\R{p_i}_*\O_{X_w},\O_{X_i})
\simeq\ext^{\ell(w)-1}_{\O_{X_i}}(\O_{p_i(X_w)},\O_{X_i})[1-\ell(w)].$$
Hence \eqref{eq:isorr} implies that
$$\R{p_i}_*\Rhom_{\O_X}(\O_{X_w},\Omega_{X/X_i}) \simeq
\ext^{\ell(w)-1}_{\O_{X_i}}(\O_{p_i(X_w)},\O_{X_i})[-\ell(w)].$$
Applying $\R{p_i}_*$ to the distinguished triangle \eqn
&&\ext^{\ell(w)}_{\O_X}(\O_{X_w},\Omega_{X/X_i})[-\ell(w)] \to
\Rhom_{\O_X}(\O_{X_w},\Omega_{X/X_i}) \to
\F\otimes\Omega_{X/X_i}\To[+1],\eneqn we obtain a distinguished triangle
\eqn&& \ba{rl}
\R{p_i}_*\ext^{\ell(w)}_{\O_X}(\O_{X_w},\Omega_{X/X_i})[-\ell(w)]
\to& \ext^{\ell(w)-1}_{\O_{X_i}}(\O_{p_i(X_w)},\O_{X_i})[-\ell(w)]\\
&\qquad\qquad\to \R{p_i}_*(\F\otimes\Omega_{X/X_i}) \To[\,+1\,]. \ea
\eneqn Hence, taking cohomology groups, we obtain
$$
R^{k}{p_i}_*(\F\otimes\Omega_{X/X_i}) \isoto
R^{k-\ell(w)+1}{p_i}_*\ext^{\ell(w)}_{\O_X}(\O_{X_w},\Omega_{X/X_i})
\simeq 0
$$
for $k>\ell(w)$.
On a neighborhood of $p_i(\dX_x)$,
$S\to X_i$ is an embedding, and we have
$R^{k}{p_i}_*(\F\otimes\Omega_{X/X_i})\simeq
{p_i}_*H^k(\F\otimes \Omega_{X/X_i})$. Hence
$H^k(\F\otimes \Omega_{X/X_i})=0$ for
$k>\ell(w)$ on a neighborhood of $\dX_x$, which is a contradiction.
\QED

\section{Affine flag manifolds}
\label{sec:aff} In this section we shall consider affine flag
manifolds. Let $(a_{ij})_{ij\in I}$ be an affine Cartan matrix. Let
$Q=\soplus_{i\in I}\Z\al_i$ and $L=\soplus_{i\in I}\Z\Lambda_i$. We
take $L\oplus Q$ as the integral weight lattice $P$. In the Cartan subalgebra
$\gt=\Hom(P,\C)$, we give the simple coroots $h_i\in\gt$ by
$$\lan h_i,\al_j\ran=a_{ij}, \text{ and }\lan h_i,\Lambda_j\ran=\delta_{ij}.$$
\begin{remark}
\bnum
\item
We have taken $L\oplus Q$ as the integral weight lattice $P$. Let $B$ denote
the associated Borel subgroup. Let $P'$ be another integral  weight lattice
satisfying \eqref{cond:weight} and let $B'$ be its associated Borel
subgroup. Then there is a map $P\to P'$ and therefore a morphism
$B'\to B$. Hence we have morphisms $K_B(\pt)\to K_{B'}(\pt)$ and
$K_{B'}(\pt)\otimes_{K_B(\pt)}K_{B}(X)\isoto K_{B'}(X)$. In this
sense, our choice of $P$ is universal.
\item
However our choice of $P$ can often be realized as a direct sum
$P'\oplus S$, where $S\subset P^W$ and $Q\subset P'$. Then we have
$K_{B}(X)\simeq \Z[S]\otimes_\Z K_{B'}(X)$. \ee
\end{remark}

The Weyl group $W$ acts on $P$ and $Q$,
and we have an exact sequence of $W$-modules
$$0\to Q\to P\to L\to0,$$
where $L$ is endowed with the trivial action of $W$.

Let us set $Q_+=\sum_{i\in I}\Z_{\ge0}\al_i$ and let
$\delta\in Q_+$ be the imaginary root such that
$$\set{\alpha\in Q}{\text{$\lan h_i,\al\ran=0$ for all $i$}}=\Z\delta.$$
Similarly, let us choose
$c\in \sum_{i\in I}\Z_{\ge0}h_i$ such that
$$\Bigl\{h\in \ssum_{i\in I}\Z h_i\,;\,
\text{$\lan h,\al_i\ran=0$ for all $i$}\Bigr\}=\Z c.$$
We write
$$\text{$\delta=\sum_{i\in I}a_i\al_i$ and $c=\sum_{i\in I}a_i^\vee h_i$.}$$
Then there exists a unique symmetric bilinear form
$(\ ,\ )\cl P\times P\to \Q$ such that
\eq&&
\ba{rl}
(\la,\la')&=0\quad \text{for any $\la,\la'\in L$,}\\[.5ex]
\lan h_i,\la\ran&=\dfrac{2(\al_i,\la)}{(\al_i,\al_i)}
\quad\text{for any $\la\in P$,}\\[2.5ex]
\lan c,\la\ran&=(\delta,\la)\quad \text{for any $\la\in P$.}
\ea
\eneq

Set $Q_\cla=Q/\Z\delta$.
Then $W$ acts on $Q_\cla$.
Let us choose $0\in I$ such that
the image $W_\cla$ of $W$ in $\Aut(Q_\cla)$ is generated by
the image of $\{s_i\}_{i\in I\setminus\{0\}}$
and $a_0=1$. These conditions are equivalent to saying that
$a_0=1$ and $\delta-\alpha_0$ is a constant multiple of a root.
Such a $0$ exists uniquely up to a Dynkin diagram automorphism.
Note that
$$(\al_0,\al_0)=
\begin{cases}
2&\text{if $\g\not\simeq A^{(2)}_{2n}$,}\\
4&\text{if $\g\simeq A^{(2)}_{2n}$,}
\end{cases}\quad
\text{and }\quad
\theta\seteq\delta-\al_0\in\dfrac{(\al_0,\al_0)}{2}\Delta^+.$$

Let $P^W$ denote the space of $W$-invariant integral weights.
We have in the affine case
$$\Z[P^W]\isoto\Z[P]^W,$$
since, for any $\la\in P$, either
$\la\in P^W$ or $W\la$ is an infinite set (see Remark~\ref{rem:inf}).

We see easily that
$P^W=\set{\la\in P}{\text{$\lan h_i,\la\ran=0$ for all $i$}}$
is isomorphic to $Q$ by the map
$\eta\cl Q\isoto P^W$ given by
\eq
\eta(\beta)=\beta-\sum_{i\in I}\lan h_i,\beta\ran\Lambda_i.
\label{eq:eta}
\eneq

\smallskip
Hence we have
$P=L\soplus P^W$ and $\Z[P]\cong\Z[L]\otimes_\Z\Z[P^W]$,
which implies
$$K_B(\pt)\otimes_{Z[P]^W}\Z[P]\cong
\Z[P]\otimes_{\Z[P^W]}\Z[P]\cong \Z[P]\otimes_\Z\Z[L]
\cong \Z[L]\otimes_\Z\Z[P].$$
For any $w\in W$,
let $i_w^*\circ\beta\cl \Z[P]\otimes_{\Z[P^W]}\Z[P]\to\Z[P]$
be as in \eqref{eq:iw}.

\Lemma\label{lem:inj}
The homomorphism induced by the $i_w^*\circ\beta$'s
$$\Z[P]\otimes_{\Z[P^W]}\Z[P]\to\Z[P]^{\prod W}$$
is injective.
Here $\Z[P]^{\prod W}$ is the product of the copies of $\Z[P]$
parameterized by the elements of $W$.
\enlemma
We shall give the proof of this lemma in \S\,\ref{sec:inj}.

As a corollary (together with \eqref{eq:embK}),
we have
\Cor\label{cor:inj}
The homomorphism
$$K_B(\pt)\otimes_{\Z[P^W]}\Z[P]\To[\beta] K_B(X)$$
is injective.
\encor

Let $R$ be the subring of $\Q(\e^\delta)$
generated by $\e^{\pm\delta}$ and $(\e^{n\delta}-1)^{-1}$ ($n>0$).
Then we have an injective homomorphism
\eq
&&\beta\cl R\otimes_{\Z[\e^{\pm\delta}]}K_B(\pt)\otimes_{\Z[P^W]}\Z[P]\mono
R\otimes_{\Z[\e^{\pm\delta}]}K_B(X).
\label{mor:R}
\eneq

Our main result is the following theorem. \Theorem\label{th:main}
For all $w\in W$, $[\O_{X_w}]\in K_B(X)$, considered as an element
of $R\otimes_{\Z[\e^{\pm\delta}]}K_B(X)$, is in the image of the map
$\beta$ in \eqref{mor:R}. \entheorem %

Hence, $\soplus_{w\in W}\hs{-1ex}K_B(\pt)[\O_{X_w}]$ may be regarded as a
submodule of
$R\hs{-.5ex}\tens_{\Z[\e^{\pm\delta}]}K_B(\pt)\otimes_{\Z[P^W]}\Z[P]$.

We shall prove a slightly more precise result.

Let $\xi\cl\Z[P]\to K_B(\pt)\otimes_{\Z[P^W]}\Z[P]\simeq\Z[L]\otimes_\Z\Z[P]$
be the homomorphism given by
\eq
\text{$\xi(\e^{\la+\beta})=\e^{-\la}\otimes \e^{\la+\beta}$
for $\la\in L$ and $\beta\in Q$.}
\label{def:xi}
\eneq%
We extend this to
$$\xi\cl
R\otimes_{\Z[\e^{\pm\delta}]}\Z[P]\mono
R\otimes_{\Z[\e^{\pm\delta}]}K_B(\pt)\otimes_{\Z[P^W]}\Z[P].
$$

Set $\rho=\sum_{i\in I}\La_i$.
Then $\kappa^*\seteq\lan c,\rho\ran=\sum_{i\in I}a_i^\vee$
is called the {\em dual Coxeter number} of $\g$.
Let us introduce the $W$-submodule of $\Z[P]$:
\eq
\Z[P]_{[0,-\kappa^*)}\seteq
\soplus_{\la\in P,\,0\ge\lan c,\la\ran>-\kappa^*}\Z\e^\la.
\label{eq:dC}
\eneq

\Prop\label{prop:main}
For any $w\in W$, there exists a unique
$$\Gro_w\in R\otimes_{\Z[\e^{\pm\delta}]}
\Z[P]_{[0,-\kappa^*)}$$
such that $\beta\circ\xi(\Gro_w)=[\O_{X_w}]$.
\enprop
We call $\Gro_w$ the {\em affine Grothendieck polynomial}.
The proof of this proposition is given in \S\,\ref{sec:main}
as an application of Theorem~\ref{th:van} below.

\begin{remark}
In the finite-dimensional case,
$Z[P]^W$ is much bigger than $\Z[P^W]$, and
the choice of Grothendieck polynomials
is not unique. However, in the affine case $\Gro_w$ is uniquely determined.
\end{remark}

\begin{remark}
Let $\psi\cl\Z[P]\to\Z[P]$ be
the homomorphism given by
$$\mbox{$\psi(\e^{\la+\al})=\e^{\la-\eta(\al)}$ for $\la\in L$ and $\al\in Q$.}
$$
Then we have
\eq \Gro_{w^{-1}}=\psi(\Gro_w)\quad\text{for any $w\in W$.} \label{eq:inv}
\eneq
Indeed, let $\phi\cl\g\to\g$ be the Lie algebra homomorphism such that
$\phi(e_i)=f_i$, $\phi(f_i)=e_i$ and $\phi(h)=-h$ for $h\in\gt$. It
induces  group scheme morphisms $\phi\cl B\to B_-$ and $\phi\cl
B_-\to B$.
Note that $\phi$ induces an isomorphism
$$\Z[P]\simeq K_T(\pt)\simeq K_B(\pt)
\isoto[{\phi^*}]K_{B_-}(\pt)\simeq K_T(\pt)
\simeq\Z[P],$$
which is given by $\e^\la\mapsto \e^{-\la}$.
There exists a unique scheme isomorphism $a\cl G\to G$ such
that $a(e)=e$ and $a(bgb_-^{-1})=\phi(b_-)a(g)\phi(b)^{-1}$ for
$b\in B$, $b_-\in B_-$ and $g\in G$.
We have $K_B(X)\simeq K_{B\times B_-}(G)$, and $a\cl G\to G$ induces
a commutative diagram
$$\xymatrix{
\Z[P]\ar[r]^(.28)\xi\ar[d]_{\psi} &K_B(\pt)\otimes
K_{B_-}(\pt)\ar[d]_k\ar[r]^(.6)\beta &K_{B\times B_-}(G)\ar[d]^{a}
\ar@{-}[r]^(.55)\sim& K_B(X)\ar[d]^{\tilde\psi}\\
\Z[P]\ar[r]^(.28)\xi &K_B(\pt)\otimes
K_{B_-}(\pt)\ar[r]^(.6)\beta&K_{B\times B_-}(G) \ar@{-}[r]^(.55)\sim
&K_B(X).}
$$
Here $k(\e^\la\otimes \e^\mu)=\e^{-\mu}\otimes \e^{-\la}$.
Indeed, $k\xi(\e^{\la+\al})=k(\e^{-\la}\otimes\e^{\la+\al})
=\e^{-\la-\al}\otimes\e^{\la}
=\e^{-\la-\al+\eta(\al)}\otimes\e^{\la-\eta(\al)}
=\xi(\e^{\la-\eta(\al)})$.

Since $a(BwB_-)=Bw^{-1}B_-$,
we have $\tilde\psi([\O_{X_w}])=[\O_{X_{w^{-1}}}]$. which implies
\eqref{eq:inv}.

\end{remark}

\section{Vanishing of the Weyl group cohomology}
\label{sec:van}

Let $(W,S)$ be a Coxeter group,
where $S$ is a system of generators.
For a subset $S'$ of $S$, let us denote by $W_{S'}$ the subgroup of $W$
generated by $S'$.
\Lemma\label{lem:van}
For any $W$-module $V$, we have
$$H^1(W,V)=\dfrac{\set{\{v_s\}_{s\in S}}{\text{$v_s\in V$ satisfy
conditions \eqref{en:com1} and \eqref{en:com2} below}}}
{\set{\{v_s\}_{s\in S}}{\text{there exists $v\in V$ such that
$v_s=(1-s)v$}}},
$$
\bnum
\item $(1+s)v_s=0$,\label{en:com1}
\item for any pair of distinct elements $s,t\in S$ such that
$W_{\{s,t\}}$ is a finite group,
$$\sum_{w\in W_{\{s,t\}},\,ws>w}(-1)^{\ell(w)}wv_s=
\sum_{w\in W_{\{s,t\}},\,wt>w}(-1)^{\ell(w)}wv_t.$$
Here $\ell\cl W\to\Z_{\ge0}$ is the length function.
\label{en:com2}
\ee
\enlemma
\Proof
Let $\Z[W]\to \Z$ be the augmentation homomorphism $W\ni w\mapsto 1$.
Then its kernel is the image of the homomorphism
$\varphi\cl \oplus_{s\in S}\Z[W]e_s\to \Z[W]$
given by $\varphi(e_s)=1-s$.
By definition, $H^1(W,V)=\Ext^1_{\Z[W]}(\Z,V)$
is the cohomology of
$$\Hom_{\Z[W]}(\Z[W],V)\to \Hom_{\Z[W]}(\oplus_{s\in S}\Z[W]e_s,V)
\to\Hom_{\Z[W]}(\Ker(\varphi),V).$$
Hence it is enough to show that
$\Ker(\varphi)$ coincides with the $\Z[W]$-submodule $N$
generated by
$(1+s)e_s$ ($s\in W$) and
$\sum_{x\in W_{\{s,t\}},\,xs>x}(-1)^{\ell(x)}xe_s
-\sum_{x\in W_{\{s,t\}},\,xt>x}(-1)^{\ell(x)}xe_t$
where $(s,t)\in S\times S$ ranges over the pairs
as in \eqref{en:com2}.
It is easy to see that those elements are in $\Ker(\varphi)$.
Indeed, the last elements belong to $\Ker(\varphi)$ because
$\sum_{x\in W_{\{s,t\}},\,xs>x}(-1)^{\ell(x)}x(1-s)
=\sum_{x\in W_{\{s,t\}}}(-1)^{\ell(x)}x$.

Let $z=\sum_{w\in W,s\in S} a_{w,s}we_s$ be an element of
$\Ker(\varphi)$  where $a_{w,s}\in\Z$.

Since $\sum_{w\in W,s\in S} a_{w,s}w(1-s)=0$, we have
\eq
&&\sum_{s\in S}a_{w,s}=\sum_{s\in S}a_{ws,s}\quad\text{for all $w\in W$.}
\label{eq:ws}
\eneq

We shall show $z\in N$ by induction on $\ell$, the largest number
among the $\ell(w)$'s such that $a_{w,s}\not=0$ for some $s\in S$.
Then we shall show $z\in N$ by induction on the cardinality of
$\set{(w,s)\in W\times S}{\text{$a_{w,s}\not=0$ and
$\ell=\ell(w)$}}$. Let us take $(w_1,s_1)$ such that
$a_{w_1,s_1}\not=0$ and $\ell=\ell(w_1)$. Subtracting
$a_{w_1,s_1}w_1(1+s_1)e_{s_1}\in N$ from $z$ when $w_1s_1<w_1$, we
may assume from the beginning that $w_2\seteq w_1s_1>w_1$. Then
$\ell(w_2)=\ell+1$. Applying \eqref{eq:ws} for $w=w_2$, we have
$0=\sum_{s\in S}a_{w_2,s}=\sum_{s\in S}a_{w_2s,s}
=a_{w_1,s_1}+\sum_{s\in S,\,s\not=s_1}a_{w_2s,s}$. Hence there
exists $s_2\not=s_1$ such that $a_{w_2s_2,s_2}\not=0$. Hence we have
$\ell(w_2s_2)\le\ell<\ell(w_2)=\ell+1$. Hence $w_2$ is the longest
element in $w_2W_{s_1,s_2}$. Therefore $W_{s_1,s_2}$ is a finite
group. Let $w_3$ be the shortest element of $w_1W_{s_1,s_2}$.
Subtracting $\pm a_{w_1,s_1}\Bigl(
\sum_{x\in W_{\{s_1,s_2\}},\,xs_1>x}(-1)^{\ell(x)}w_3xe_{s_1}
-\sum_{x\in W_{\{s_1,s_2\}},\,xs_2>x}(-1)^{\ell(x)}w_3xe_{s_2}\Bigr)$
in $N$ from $z$,
we can erase the term $w_1e_{s_1}$ in $z$, and
the induction proceeds.
\QED

Now let us return to the affine case
where $W$ is the Weyl group.
Recall that $R$ is the ring generated by $\e^{\pm\delta}$
and $(\e^{n\delta}-1)^{-1}$ ($n\not=0$).

Note that $\Z[P]$ is a direct sum of $W$-submodules of the form
$\soplus_{\la\in W\la_0}\Z\e^\la$ ($\la_0\in P$).

\Theorem\label{th:van}
\bnum
\item
If $\vert I\vert>2$, then
$H^1(W,\,R\otimes_{\Z[\e^{\pm\delta}]}\Z[P])=0$.\label{th:van1}
\item
For any affine Lie algebra $\g$,
$H^1\bigl(W,R\otimes_{\Z[\e^{\pm\delta}]}
(\hs{-1ex}\soplus_{\vert\lan c,\la\ran\vert<\lan c,\rho\ran}\Z\e^\la)\bigr)=0$
where $\rho=\sum\limits_{i\in I}\Lambda_i$.\label{th:van2}
\ee
\entheorem

In fact, we shall prove more precise results.
For $J\subset I$, let $W_J$ be the subgroup of $W$
generated by $\set{s_i}{i\in J}$.

\Prop\label{prop:van}
\bnum
\item
If $J\not=I$, then $H^1(W_J,\Z[P])=0$,\label{en:van1}
\item
if $\lan c,\la_0\ran\not=0$
and if $\la_0$ satisfies one of the conditions below
\be[{\rm(a)}]
\item
$\la_0$ is not regular
\ro{\em i.e.}, $(\beta,\la_0)=0$ for some $\beta\in\Delta$\rf,
\item $\vert I\vert>2$,
\ee
then
$H^1(W,\,\oplus_{\la\in W\la_0}\Z\e^\la)=0$,\label{en:van2}
\item
if $\lan c,\la_0\ran=0$,
and $\lan h_i,\la_0\ran\ge0$ for all $i\in I\setminus\{0\}$, then
$$1-\e^{\lan h_0,\la_0\ran\delta}\cl
H^1(W,\,\oplus_{\la\in W\la_0}\Z\e^\la)
\to H^1(W,\,\oplus_{\la\in W\la_0}\Z\e^\la)
$$
is the zero map.\label{en:van3} \ee
\enprop %
Note that in \eqref{en:van3}, we have $W\la_0+\Z\lan
h_0,\la_0\ran\delta=W\la_0$. Also note that, if $0<\vert\lan
c,\la_0\ran\vert<\lan c,\rho\ran$, then $\la_0$ is not regular.
Together with $H^1(W,\Z)=0$, it is easy to see that
Proposition~\ref{prop:van} implies Theorem~\ref{th:van}. The proof
of Proposition~\ref{prop:van} will be given in the next section.

\begin{remark}
If $\vert I\vert=2$ and $\la_0$ is regular dominant, then
$H^1(W,\soplus_{\la\in W\la_0}\Z\e^\la)\simeq\Z$.
\end{remark}

\section{Proof of Proposition \ref{prop:van}}
\label{sec:prvan}

In this section, we shall prove Proposition \ref{prop:van}.

\medskip\noindent
\subsection{Proof of Proposition \ref{prop:van} (i), (ii)}
In the case \eqref{en:van2}, we may assume that $\lan c,\la_0\ran>0$.
In such a case, $W\la_0$ contains a dominant weight.
Hence, in order to prove \eqref{en:van1} and \eqref{en:van2},
it is enough to show that
\eq
&&
H^1(W_J,\soplus_{\la\in W_J\la_0}\Z\e^\la)=0
\label{st:vanh}
\eneq
under the condition
\eq
&&\parbox{65ex}
{if $J\subset I$ and $\lan h_i,\la_0\ran\ge0$ for any $i\in J$.
Moreover, when $J=I$ and $\vert I\vert=2$, we assume further that
$\la_0$ is not regular.}\label{cond:Jlam}
\eneq

We shall show this by induction on the cardinality of $J$.
If $\vert J\vert\le1$, then it is obvious.
Assuming that $\vert J\vert>1$, let us take $i_0\in J$,
and set $J_0=J\setminus\{i_0\}$. Then \eqref{st:vanh}
is true for $J_0$ by the induction hypothesis.

Assuming that $(v_i)_{i\in J}$ with $v_i\in \soplus_{\la\in
W_J\la_0}\Z\e^\la$ satisfies \eqref{en:com1} and \eqref{en:com2} in
Lemma \ref{lem:van}, let us show the existence of $v\in
\soplus_{\la\in W_J\la_0}\Z\e^\la$ such that $v_i=(1-s_i)v$. By the
induction hypothesis, there exists $v'\in \soplus_{\la\in
W_J\la_0}\Z\e^\la$ such that $v_i=(1-s_i)v'$ for all $i\in J_0$.
Hence replacing $v_i$ with $v_i-(1-s_i)v'$, we may assume from the
beginning that $v_i=0$ for all $i\in J_0$. On the other hand, since
$(1+s_{i_0})v_{i_0}=0$, there exists $u_0$ such that
$v_{i_0}=(1-s_{i_0})u_0$. By \eqref{en:com2} in Lemma \ref{lem:van},
we have \eq\label{eq:u0} &&\text{$\sum_{w\in
W_{\{i_0,j\}}}(-1)^{\ell(w)}wu_0=0$ if $j\in J_0$ and
$W_{\{i_0,j\}}$ is a finite group.} \eneq It is enough to show the
following: \eq\label{eq:z0z1} &&\ba{l} \text{there exists a
decomposition
 $u_0=z_0+z_1$}\\
\text{where
$z_0,z_1\in \soplus_{\la\in W_J\la_0}\Z\e^\la$
and $s_{i_0}z_0=z_0$ and $z_1$ is $W_{J_0}$-invariant.}
\ea
\eneq
Indeed, we  then have
$v_{i_0}=(1-s_{i_0})u_0=(1-s_{i_0})z_1$
and $v_i=0=(1-s_i)z_1$ for $i\in J_0$.

\medskip
We shall show \eqref{eq:z0z1} under the conditions \eqref{cond:Jlam}
and \eqref{eq:u0}. Let $d\cl W_J\la_0\to \Z_{\ge0}$ be the function
given by $d(\la)=\sum_{i\in J}m_i$ writing $\la_0-\la=\sum_{i\in
J}m_i\al_i$ ($m_i\in\Z_{\ge0}$). Let us write $u_0=\sum_{\la\in
W_J\la_0}a_\la\e^\la$. Set $\supp(u_0)\seteq\set{\la\in
W_J\la_0}{a_\la\not=0}$, and we argue by induction on
$d(u_0)\seteq\max\set{d(\la)}{\la\in\supp(u_0)}$. Then we argue by
induction on the cardinality of
$\supp^{\max}(u_0)\seteq\set{\la\in\supp(u_0)}{d(\la)=d(u_0)}$.

Let us take $\la_1\in\supp^{\max}(u_0)$.
If $\lan h_{i_0},\la_1\ran< 0$,
then $d(s_{i_0}\la_1)<d(\la_1)$ and, subtracting the $s_{i_0}$-invariant
$a_{\la_1}(\e^{\la_1}+\e^{s_{i_0}\la_1})$ from $u_0$,
we can delete the term $\e^{\la_1}$ in $u_0$.
If $\lan h_{i_0},\la_1\ran= 0$, then
subtracting the $s_{i_0}$-invariant
$a_{\la_1}\e^{\la_1}$ from $u_0$,
we can delete the term $\e^{\la_1}$ in $u_0$.
Hence we can assume that $\lan h_{i_0},\la_1\ran>0$.

Now assume that $\lan h_i,\la_1\ran>0$ for some $i\in J_0$. Then
$\la_1$ is regular dominant with respect to $\{i_0,i\}$. If
$I=\{i_0,i\}$, then $\la_1=\la_0$ and it contradicts the hypothesis
that $\la_0$ is not regular. Hence $I\not=\{i_0,i\}$ and
$W_{\{i_0,i\}}$ is a finite group. Therefore \eqref{eq:u0} implies
that $\sum_{w\in W_{\{i_0,i\}}}(-1)^{\ell(w)}a_{w\la_1}=0$. For
$w\in W_{\{i_0,i\}}\setminus\{e\}$, we have $a_{w\la_1}=0$, because
$d(w\la_1)>d(\la_1)$. Thus we obtain the contradiction
$a_{\la_1}=0$.

We thus conclude $\lan h_i,\la_1\ran\le0$ for all $i\in J_0$.
Hence $d(\la)<d(\la_1)$ for all $\la\in W_{J_0}\la_1\setminus\{\la_1\}$.
Then subtracting the $W_{J_0}$-invariant
$a_{\la_1}\bigl(\suml_{\la\in W_{J_0}\la_1}\e^{\la}\bigr)$
from $u_0$, we can erase the term $\e^{\la_1}$ in $u_0$, and
the induction proceeds.
Note that $W_{J_0}$ is a finite group.

Thus we have proved \eqref{eq:z0z1} under the conditions
\eqref{cond:Jlam} and \eqref{eq:u0}.

\bigskip\noindent
\subsection{Proof of Proposition \ref{prop:van} (iii)} 
We may assume that $\lan h_0,\la_0\ran<0$.
Set $I_0=I\setminus\{0\}$,
and let $W_0$ be the subgroup of $W$ generated by $\{s_i\}_{i\not=0}$.
Then $W\la_0\subset W_0\la_0+\Z\delta$.
Hence,
as in the proof of \eqref{en:van2}, it is enough to show that
if $u_0\in\soplus_{\la\in W_0\la_0}\Z[\e^{\pm\delta}]
\e^\la$ satisfies the condition
\eq\label{eq:u00}
&&\text{$\sum_{w\in W_{\{0,j\}}}(-1)^{\ell(w)}wu_0=0$
if $j\in I_0$ and $W_{\{0,j\}}$ is a finite group,}
\eneq
then
\eq\label{eq:z0z10}
&&\ba{l}
\text{there exists a decomposition
 $(1-\e^{\lan h_0,\la_0\ran\delta})u_0=z_0+z_1$}\\
\text{where
$z_0,z_1\in \soplus_{\la\in W_0\la_0}\Z[\e^{\pm\delta}]\e^\la$
and $s_{0}z_0=z_0$ and $z_1$ is $W_0$-invariant.}
\ea
\eneq
Let us write
$u_0=\sum_{\la\in W_0\la_0}a_\la\e^{\la}$
with $a_\la\in\Z[\e^{\pm\delta}]$.
Let us set $\theta=\delta-\al_0\in \frac{(\al_0,\al_0)}{2}\Delta^+$.
Let $s_\theta$ be the reflection with respect to
$\theta$\;: $s_\theta(\la)=\la-\frac{2(\theta,\la)}{(\theta,\theta)}\theta$.
Then $s_\theta$ belongs to $W_0$, and
$s_0\la_0=s_\theta\la_0-\lan h_0,\la_0\ran\delta$.
We have an $s_0$-invariant
$z_0\seteq \e^{\la_0}+\e^{s_0\la_0}=\e^{\la_0}
+\e^{-\lan h_0,\la_0\ran\delta}\e^{s_\theta\la_0}$.
Subtracting a constant multiple of $z_0$ from $u_0$,
we may assume that $a_{s_\theta\la_0}$ vanishes.
On the other hand,
$z_1\seteq\sum_{\la\in W_0\la_0}\e^\la$ is a $W_0$-invariant.
Their linear combination $\e^{-\lan h_0,\la_0\ran\delta}z_1-z_0$
has no term $\e^{s_\theta\la_0}$ and the coefficient of
$\e^{\la_0}$ is $\e^{-\lan h_0,\la_0\ran\delta}-1$.
Hence subtracting a constant multiple of it from
$(1-\e^{\lan h_0,\la_0\ran\delta})u_0$,
we may assume that $a_{\la_0}$ and $a_{s_\theta\la_0}$ vanish.
Let us set
$\supp(u_0)=\set{\la\in W_0\la_0}{a_\la\not=0}$.
By subtracting an $s_0$-invariant from $u_0$,
we may assume further that
\eq&&
\text{$\lan \theta,\la\ran=-\lan h_0,\la\ran>0$ for any $\la\in\supp(u_0)$.}
\label{eq:u0s0}
\eneq

Hence we have reduced the problem to proving
\eq
&&\text{if $a_{\la_0}=a_{s_\theta\la_0}=0$
and if $u_0$ satisfies \eqref{eq:u00} and \eqref{eq:u0s0},
then $u_0=0$.}
\eneq

If $\vert I\vert$ is $2$, it is obvious, since
$W_0\la_0=\{\la_0,s_\theta\la_0\}$.
Let us assume $\vert I\vert>2$. Hence $W_{\{0,i\}}$
is a finite group for all $i\in I_0$.

For $\la\in W_0\la_0$, we set
\eq
a_{\la+n\delta}=\e^{-n\delta}a_\la,\label{eq:adelta}
\eneq
so that we have
$a_{\la+n\delta}\e^{\la+n\delta}=a_{\la}\e^{\la}$.
Then \eqref{eq:u00} reads as
\eq
&&\text{$\sum_{w\in W_{\{0,i\}}}(-1)^{\ell(w)}a_{w\la}=0$
for any $\la\in W_0\la_0+\Z\delta$ and $i\in I_0$.}
\label{eq:u0w0}
\eneq
Note that $s_0\la=s_\theta\la-\lan h_0,\la\ran\delta$ and
\eq
a_{s_0\la}=\e^{\lan h_0,\la\ran\delta}a_{s_\theta\la}.
\eneq

\Sub\label{sub1}
Let $\la\in W_0\la_0$.
If $k\in I_0$ satisfies $(\al_0,\al_k)=0$,
then $a_\la=a_{s_k\la}$.
\ensub
\Proof
We may assume that $\lan h_0,\la\ran<0$.
Then, \eqref{eq:u0s0} implies that
$W_{\{0,k\}}\la\cap\bigl(\supp(u_0)+\Z\delta\bigr)
\subset\{\la,s_k\la\}$,
and \eqref{eq:u0w0} implies the desired result.
\QED
Set
$$I_1\seteq\set{k\in I}{(\al_0,\al_k)=0},$$
and let $W_1$ be the subgroup of $W$ generated by
$\set{s_k}{k\in I_1}$.
Then Sublemma~\ref{sub1} implies that
\eq
&&\text{$a_\la=a_{w\la}$ for any $w\in W_1$.}\label{eq:W1inv}
\eneq

Now we shall divide the proof into two cases:
\be[{\rm (A)}]
\item there exists a $1\in I_0\setminus I_1$ such that
$(\al_1,\al_1)\not=(\al_0,\al_0)$,\label{caseA}
\item
for all $i\in I_0\setminus I_1$,
we have $(\al_i,\al_i)=(\al_0,\al_0)$.\label{caseB}
\ee

\medskip\noindent
{\em Case\eqref{caseA}}\quad
In this case, $I=\{1\}\sqcup I_1$ and
$\lan h_0,\al_1\ran\lan h_1,\al_0\ran=2$
as seen by the classification of affine Dynkin diagrams.
Note that $\{0,1\}$ is a Dynkin diagram of type {C$_2$}.
Then $\theta\seteq\delta-\al_0=\sum_{i\not=0}a_i\al_i$
satisfies $\lan h_0,\theta\ran=-2$, and hence we have
$\lan h_0,\al_1\ran a_1=-2$, which implies
$a_1+\lan h_1,\al_0\ran=0$.
Since $\lan h_1,\theta\ran=-\lan h_1,\al_0\ran$,
we have $\beta\seteq s_1\theta
=(a_1+\lan h_1,\al_0\ran)\al_1+\sum_{i\in I_1}a_i\al_i
=\sum_{i\in I_1}a_i\al_i$.
Hence $\beta$ is a constant multiple of a root
in $\Delta\cap(\sum_{i\in I_1}\Z\al_i)$, and $s_\beta$ belongs to $W_1$.
Assuming that $u_0$ does not vanish, let us choose
an element $\mu$ in $\supp(u_0)$,
highest with respect to $I_0$ ({\em i.e.}, maximal with respect to
the ordering $\ge$:
$\mu\ge\mu'$ if $\mu-\mu'\in\sum_{i\in I_0}\Z_{\ge0}\al_i$).
By \eqref{eq:W1inv},
we have
\eq
&&\text{$\lan h_k,\mu\ran\ge 0$ for any $k\in I_1$.}\label{eq:mudom}
\eneq
Let us show that $\lan h_1,\mu\ran\ge 0$.
Otherwise, $\mu$ is regular and anti-dominant with respect to $\{0,1\}$.
By \eqref{eq:u0s0}, we have
$W_{\{0,1\}}\mu\cap(\supp(u_0)+\Z\delta)
\subset\{\mu,s_1\mu,s_1s_0\mu,s_1s_0s_1\mu\}$.
We have $s_1\mu>\mu$, and hence $a_{s_1\mu}=0$.
Since $s_\beta=s_1s_\theta s_1$, we have
$a_{s_1s_\theta\mu}=a_{s_\beta s_1\mu}=a_{s_1\mu}=0$
by \eqref{eq:W1inv}.
Hence we have $a_{s_1s_0\mu}=0$.
Thus we obtain $W_{\{0,1\}}\mu\cap(\supp(u_0)+\Z\delta)
\subset\{\mu,,s_1s_0s_1\mu\}$.
Hence \eqref{eq:u0w0} implies that
$a_\mu-a_{s_1s_0s_1\mu}=0$.
On the other hand, we have
$a_{s_1s_0s_1\mu}=\e^{\lan h_0,s_1\mu\ran\delta}a_{s_1s_\theta s_1\mu}
=\e^{\lan h_0,s_1\mu\ran\delta}a_{s_\beta\mu}
=\e^{\lan h_0,s_1\mu\ran\delta}a_{\mu}$.
Hence $(1-\e^{\lan h_0,s_1\mu\ran\delta})a_\mu=0$.
Since $\lan h_0,s_1\mu\ran=\lan s_1h_0,\mu\ran<0$, we obtain
$a_\mu=0$, which is a contradiction.

We thus conclude that $\lan h_1,\mu\ran\ge0$. Along with \eqref{eq:mudom},
$\mu\in W_0\la_0$ is dominant with respect to $I_0$, and hence
we conclude $\mu=\la_0$, which contradicts $a_{\la_0}=0$.

\medskip\noindent
{\em Case\eqref{caseB}}\quad
The proof in this case is similar to the one in Case\eqref{caseB},
but slightly more complicated.
In this case, $\vert I_0\setminus I_1\vert$
is one or two by the classification of affine Dynkin diagrams.
The case $\vert I_0\setminus I_1\vert=2$ is exactly the case $A^{(1)}_n$
($n\ge2$).
Set $I_0\setminus I_1=\{i_1,i_2\}$
(when $\vert I_0\setminus I_1\vert=1$, by convention $i_1=i_2$).
We have $(\al_0,\al_i)=-1$ for $i\in I_0\setminus I_1$,
$\theta\seteq\delta-\al_0$ is a root,
and
$$\delta=\al_0+\al_{i_1}+\al_{i_2}+\sum_{i\in I_1}a_i\al_i,
\quad\theta=\al_{i_1}+\al_{i_2}+\sum_{i\in I_1}a_i\al_i$$
(when $\vert I_0\setminus I_1\vert=1$,
$\delta=\al_0+2\al_{i_1}+\sum_{i\in I_1}a_i\al_i$).

Let $w\in W_1$ be the longest element of $W_1$.

\Sub
We have $ws_{i_1}\theta=\al_{i_2}$.
\ensub
\Proof
We have $s_{i_1}\theta=\theta-\al_{i_1}=\al_{i_2}+\sum_{i\in I_1}a_i\al_i$.
Moreover, we have $\lan h_k,s_{i_1}\theta\ran
=-\bigl\lan h_k-\lan h_k,\al_{i_1}\ran h_{i_1},\al_0\bigr\ran\ge0$
for $k\in I_1$.
Hence $s_{i_1}\theta$ is dominant with respect to $I_1$.
Hence $ws_{i_1}\theta$ is anti-dominant with respect to $I_1$:
$\lan h_k,ws_{i_1}\theta\ran\le0$ for any $k\in I_1$.
Write
$$ws_{i_1}\theta=\al_{i_2}+\beta.$$
Since $ws_{i_1}\theta$ is a root, $\beta$ has the form
$$\text{$\beta=\sum_{i\in I_1}m_i\al_i$ with $m_i\in\Z_{\ge0}$.}$$
Hence we have
$(\beta,ws_{i_1}\theta)\le0$.
On the other hand,
$(\al_0,\al_0)=(\al_{i_2}+\beta,\al_{i_2}+\beta)$ implies that
$(\beta,\beta)+2(\beta,\al_{i_2})=0$.
Hence we obtain
$(\beta,\beta)=2(\beta,\al_{i_2}+\beta)\le0$,
which implies $\beta=0$.
\QED
As a corollary, we have
$ws_{i_1}s_\theta s_{i_1}w=s_{i_2}$ and
\eq
&&s_{i_2}ws_{i_1}s_\theta=ws_{i_1}.\label{eq:wthet}\\
&&s_{i_2}s_\theta ws_{i_1}s_\theta=w\label{eq:w2}
\eneq
Indeed, the last equality follows from
$(ws_{i_2}s_\theta)(ws_{i_1}s_\theta)=
(s_{i_1}ws_{i_2})(s_{i_2}ws_{i_1})=e$.

Assuming that $u_0$ does not vanish, let us choose
$\mu$ in $\supp(u_0)$,
highest with respect to $I_0$.
By \eqref{eq:W1inv},
we have
\eq
&&\text{$\lan h_k,\mu\ran\ge 0$ for any $k\in I_1$.}\label{eq:betadom}
\eneq

Let us show that $\lan h_{i_1},\mu\ran\ge 0$.
Assume the contrary: $\lan h_{i_1},\mu\ran<0$.
Then $\mu$ is regular anti-dominant
with respect to $\{0,i_1\}$.
Since $s_{i_1}\mu>\mu$,
we have
\eq
a_{s_{i_1}\mu}=0.\label{eq:simu}
\eneq
The property \eqref{eq:u0s0} implies that
$W_{\{0,{i_1}\}}\mu\cap(\supp(u_0)+\Z\delta)
\subset\{\mu,s_{i_1}\mu,s_{i_1}s_0\mu\}$,
and hence by \eqref{eq:u0w0}, together with \eqref{eq:simu},
we have
\eq
a_{\mu}+a_{s_{i_1}s_0\mu}=0.
\eneq
Set $\mu_1\seteq ws_{i_1}s_\theta\mu$.
Then we have
$(\al_0,\mu_1)=(\al_0,s_{i_1}s_0\mu)=(s_0s_{i_1}\al_0,\mu)=(\al_{i_1},\mu)<0$
and $(\al_{i_2},\mu_1)=( s_{i_1}w\al_{i_2},s_\theta\mu)
=(\theta,s_\theta\mu)=-(\theta,\mu)<0$.
Hence $\mu_1$ is also regular anti-dominant
with respect to $\{0,i_2\}$.
By \eqref{eq:u0s0}, we have
$W_{\{0,i_2\}}\mu_1\cap(\supp(u_0)+\Z\delta)
\subset\{\mu_1,s_{i_2}\mu_1,s_{i_2}s_0\mu\}$.
Since $s_{i_2}\mu_1=s_{i_2}ws_{i_1}s_\theta\mu=ws_{i_1}\mu$
by \eqref{eq:wthet},
we have
$a_{s_{i_2}\mu_1}=a_{ws_{i_1}\mu}=a_{s_{i_1}\mu}=0$ by \eqref{eq:simu}.
Here we used \eqref{eq:W1inv} in the second equality.
Hence we have
$W_{\{0,i_2\}}\mu_1\cap(\supp(u_0)+\Z\delta)
\subset\{\mu_1,s_{i_2}s_0\mu_1\}$,
and \eqref{eq:u0w0} implies that
\eqn
&&a_{\mu_1}+a_{s_{i_2}s_0\mu_1}=0.\label{eq:mus1mu}
\eneqn
By \eqref{eq:w2}, we have $s_{i_2}s_\theta\mu_1
=s_{i_2}s_\theta ws_{i_1}s_\theta\mu=w\mu$,
which implies $a_{s_{i_2}s_\theta\mu_1}=a_\mu$ by \eqref{eq:W1inv}.
By \eqref{eq:adelta}, we have
$a_{s_{i_2}s_0\mu_1}=\e^{\lan h_0,\mu_1\ran}a_{s_{i_2}s_\theta\mu_1}
=\e^{\lan h_{i_1},\mu\ran}a_{\mu}$,
and
$a_{\mu_1}=a_{s_{i_1}s_\theta\mu}=\e^{-\lan h_0,\mu\ran}a_{s_{i_1}s_0\mu}$.
Thus we obtain
$\e^{\lan h_{i_1},\mu\ran}a_{\mu}+\e^{-\lan h_{0},\mu\ran}a_{s_{i_1}s_0\mu}=0$.
Together with \eqref{eq:mus1mu}
and $\lan h_{0},\mu\ran+\lan h_{i_1},\mu\ran<0$,
we conclude that $a_\mu=0$.
It is a contradiction.

Hence we have obtained $\lan h_{i_1},\mu\ran\ge0$.
Similarly we have $\lan h_{i_2},\mu\ran\ge0$.
Thus $\mu$ is dominant with respect to $I_0$,
and hence $\mu=\la_0$, which is a contradiction.

\section{Proof of Proposition~\ref{prop:main}}
\label{sec:main}

In this section we shall prove Proposition~\ref{prop:main} as an
application of Theorem~\ref{th:van}.
Corollary~\ref{cor:dem} implies that for any $w\in W$ we have
\eq
&&p_i^*{p_i}_*([\O_{X_w}])=
\begin{cases}
[\O_{X_{ws_i}}]&\text{if $ws_i<w$,}\\[1ex]
[\O_{X_{w}}]&\text{if $ws_i>w$.}
\end{cases}
\label{eq:demO} \eneq \Lemma \label{lem:uniq} Let $J\subset I$. If
$A\in R\otimes_{\Z[\e^{\pm\delta}]}K_B(X)$ satisfies the conditions:
\bnum
\item
$p_i^*{p_i}_*A=0$ for all $i\in J$,
\item
$p_i^*{p_i}_*A=A$ for all $i\in I\setminus J$,
\item
$i_e^*A=0$, \ee then $A=0$. \enlemma \Proof Write $A=\sum_{w\in
W}a_w[\O_{X_w}]$ (infinite sum). Then the condition (iii) implies
that $a_e=0$. Let us show $a_w=0$ by induction on $\ell(w)$.

By \eqref{eq:demO}, we have
$$p_i^*{p_i}_*A=\sum_{x\in W,\,xs_i>x}(a_x+a_{xs_i})[\O_{X_x}].$$
If $xs_i<x$ for some $i\in I\setminus J$, then (ii) implies $a_x=0$.
The condition (i) implies that $a_x+a_{xs_i}=0$ for any $x$ and $i\in J$.

Let us take $i$ such that $ws_i<w$.
If $i\notin J$, then $a_w=0$.
If $i\in J$, then the induction hypothesis implies that
$a_w=-a_{ws_i}$ vanishes.
\QED

Let us recall that we have a monomorphism (Corollary~\ref{cor:inj}):
$$
\beta\cl R\otimes_{\Z[\e^{\pm\delta}]}K_B(\pt)\otimes_{\Z[P^W]}\Z[P]\mono
R\otimes_{\Z[\e^{\pm\delta}]}K_B(X).
$$
Let us also recall
$1\otimes\sD_i\in\End(R\otimes_{\Z[\e^{\pm\delta}]}K_B(\pt)\otimes_{\Z[P^W]}\Z[P])$
which acts on the last factor $\Z[P]$ as in \eqref{eq:demz}.
Then, we have a commutative diagram (for $\xi$ see \eqref{def:xi}):
$$\xymatrix{
R\otimes_{\Z[\e^{\pm\delta}]}\Z[P]\ \ar@{>->}[r]^(.35){\xi}\ar[d]_{\sD_i}&
R\otimes_{\Z[\e^{\pm\delta}]}K_B(\pt)\otimes_{\Z[P^W]}\Z[P]\
\ar@{>->}[r]^(.6)\beta\ar[d]_{{1\otimes\sD_i}}&
R\otimes_{\Z[\e^{\pm\delta}]}K_B(X)\ar[d]_{{p_i^*{p_i}_*}}\\
R\otimes_{\Z[\e^{\pm\delta}]}\Z[P]\ \ar@{>->}[r]^(.35){\xi}&
R\otimes_{\Z[\e^{\pm\delta}]}K_B(\pt)\otimes_{\Z[P^W]}\Z[P]\
\ar@{>->}[r]^(.6)\beta&
R\otimes_{\Z[\e^{\pm\delta}]}K_B(X).
}$$
Let $j_w\cl R\otimes_{\Z[\e^{\pm\delta}]}\Z[P]\to
R\otimes_{\Z[\e^{\pm\delta}]}\Z[Q]$
be the homomorphism given by
\eq
j_w(\e^{\la+\al})=\e^{w(\la+\al)-\la}
\quad\text{for $\la\in L$ and $\al\in Q$.}\label{def:jw}
\eneq
Then we have a commutative diagram
\eq&&
\ba{c}\xymatrix{
R\otimes_{\Z[\e^{\pm\delta}]}\Z[P]\ \ar[d]^{j_w}\ar@{>->}[r]^(.35)\xi
&R\otimes_{\Z[\e^{\pm\delta}]}K_B(\pt)\otimes_{\Z[P^W]}\Z[P]
\ar[d]_{{i_w^*\circ\beta}}\\
R\otimes_{\Z[\e^{\pm\delta}]}\Z[Q]\ \ar@{^{(}->}[r]
&R\otimes_{\Z[\e^{\pm\delta}]}\Z[P].
}\ea\label{com:jw}\eneq
We have
$i_e^*[\O_{X_w}]=0$ for $w\not=e$.

Hence in order to prove Proposition~\ref{prop:main},
it is enough to construct
$\Gro_w\in R\otimes_{\Z[\e^{\pm\delta}]}\Z[P]_{[0,-\kappa^*)}$
which satisfies
\eq\label{eq:Aw}
&&\left\{\ba{lll}
{\rm(i)}\quad&\Gro_e=1,\\
{\rm(ii)}& \sD_i(\Gro_w)=
\begin{cases}
\Gro_{ws_i}&\text{if $ws_i<w$,}\\
\Gro_w&\text{if $ws_i>w$,}
\end{cases}
\quad&\text{for $w\not=e$,}\\[3ex]
{\rm(iii)}&j_e(\Gro_w)=0&\text{for $w\not=e$.}
\ea\right.
\eneq
Then Lemma~\ref{lem:uniq} guarantees that
$$\beta\circ\xi(\Gro_w)=[\O_{X_w}].$$
Hence $\Gro_w$ is the affine Grothendieck polynomial.


\medskip
We shall construct
such $\Gro_w$'s by induction on $\ell(w)$.
Assuming that
$\Gro_x\in R\otimes_{\Z[\e^{\pm\delta}]}
\Z[P]_{[0,-\kappa^*)}$ has been constructed for $x<w$ satisfying
\eqref{eq:Aw},
let us construct
$\Gro_w$.
Note that $\Gro_x$ is $s_i$-invariant if $xs_i>x$, $x<w$,
and $\sD_i\Gro_x=\Gro_{xs_i}$ if $xs_i<x<w$.

Let us set $J=\set{i\in I}{ws_i<w}$, and
$\rho_J\seteq\sum_{i\in J}\Lambda_i$.
Set $B=\e^{\rho_J}\Gro_w$.
We have
\eq
&&
\ba{l}
\text{$\sD_i\circ\e^{-\rho_J}=\e^{-\rho_J}(1-\e^{-\al_i})^{-1}\circ(1-s_i)$%
\quad for $i\in J$,}\\[1ex]
\text{$(\sD_i-1)\circ\e^{-\rho_J}
=\e^{-\rho_J-\al_i}(1-\e^{-\al_i})^{-1}\circ(1-s_i)$\quad for $i\in I\setminus J$.}
\ea
\eneq
Hence the condition \eqref{eq:Aw} (ii) reads as
\eq
(1-s_i)B=
\begin{cases}
\e^{\rho_J}(1-\e^{-\al_i})\Gro_{ws_i}&\text{if $i\in J$,}\\
0&\text{otherwise.}
\end{cases}
\label{eq:B}
\eneq
Assume that this equation is solved with
\eq &&
B\in R\tens_{\Z[\e^{\pm\delta}]}
\bigl(\oplusl_{\lan c,\rho_J\ran\ge\lan c,\la\ran
>\lan c,\rho_J-\rho\ran}\Z\e^\la\bigr).\label{eq:Bcox}
\eneq
Set $C=j_e(B)\in
R\otimes_{\Z[\e^{\pm\delta}]}\Z[Q]$.
Then
$\eta(C)$ (see \eqref{eq:eta}) belongs to
$R\otimes_{\Z[\e^{\pm\delta}]}\Z[P^W]$
and satisfies $j_e(\eta(C))=C$.
Therefore, $\Gro_w\seteq \e^{-\rho_J}(B-\eta(C))$
belongs to $R\tens_{\Z[\e^{\pm\delta}]}
\bigl(\oplusl_{0\ge\lan c,\la\ran
>-\lan c,\rho\ran}\Z\e^\la\bigr)$ and
satisfies all the conditions in \eqref{eq:Aw}.

Thus we reduced the problem to solving the equation
\eqref{eq:B} with \eqref{eq:Bcox}.

\medskip
In order to solve
\eqref{eq:B} with \eqref{eq:Bcox}, let us apply Theorem~\ref{th:van}
\eqref{th:van2}.
Note that $\e^{\rho_J}(1-\e^{-\al_i})\Gro_{ws_i}
\in R\tens_{\Z[\e^{\pm\delta}]}
\bigl(\oplusl_{\lan c,\rho_J\ran\ge\lan c,\la\ran
>\lan c,\rho_J-\rho\ran}\Z\e^\la\bigr)$.
Since $J\not=I$, we have
$\lan c,\rho\ran>\lan c,\rho_J\ran$ and $\lan c,\rho_J-\rho\ran\ge
-\lan c,\rho\ran$.
Hence, by Theorem~\ref{th:van} \eqref{th:van2},
it is enough to show the compatibility conditions
\eqref{en:com1} and \eqref{en:com2} in Lemma~\ref{lem:van}:
namely
\eq
\label{eq:last}
&&\ \left\{\parbox{30em}{
\bnum
\item
$(1+s_i)\e^{\rho_J}(1-\e^{-\alpha_i})\Gro_{ws_i}=0$ for all $i\in J$,

\vspace{.5ex}
\item
if $i,j\in J$, then
$$\hs{-4ex}\sum_{x\in W_{\{i,j\}},\,xs_i>x}\hs{-4ex}(-1)^{\ell(x)}x
\e^{\rho_J}(1-\e^{-\alpha_i})\Gro_{ws_i}
=\hs{-3.5ex}\sum_{x\in W_{\{i,j\}},\,xs_j>x}\hs{-4ex}(-1)^{\ell(x)}x
\e^{\rho_J}(1-\e^{-\alpha_j})\Gro_{ws_j},$$
\item
if $i\in J$, $j\in I\setminus J$ and $W_{\{i,j\}}$ is a finite group, then
$$\sum_{x\in W_{\{i,j\}},\,xs_i>x}\hs{-2ex}(-1)^{\ell(x)}x
\e^{\rho_J}(1-\e^{-\alpha_i})\Gro_{ws_i}=0.$$
\ee
}\right.
\eneq

\medskip\noindent
{\em Proof of \eqref{eq:last}} (i)\quad
This follows from the fact that $\Gro_{ws_i}$ is $s_i$-invariant,
which implies that $\e^{\rho_J}(1-\e^{-\alpha_i})\Gro_{ws_i}=
(1-s_i)(\e^{\rho_J}\Gro_{ws_i})$.

\medskip
In order to prove  \eqref{eq:last} (ii), (iii), let us recall the
following well-known results on Demazure operators. For $x\in W$,
$\sD_{i_1}\cdots\sD_{i_m}$ does not depend on the choice of reduced
expressions $x=s_{i_1}\cdots s_{i_m}$. We denote it by $\sD_x$. We
have $\sD_x\sD_i=\sD_x$ if $xs_i<x$.

\medskip\noindent
{\em Proof of \eqref{eq:last}} (ii)
\quad
Since $ws_i<w$ and $ws_j<w$, $w$ is the longest element in $wW_{\{i,j\}}$.
Let $w_1$ be the shortest element of
$wW_{\{i,j\}}$, and $w_0$ the longest element of $W_{\{i,j\}}$.
Hence $w=w_1w_0$.
Let $\Delta^+_{\{i,j\}}=\Delta^+\cap(\Z\al_i+\Z\al_j)$,
and $D=\prod_{\al\in \Delta^+_{\{i,j\}}}(1-\e^{-\al})$.
Then we have (see e.g.\ \cite{D,Kumar})
$$\sD_{w_0}=\e^{-\rho_J}D^{-1}\sum_{x\in W_{\{i,j\}}}(-1)^{\ell(x)}x\circ
\e^{\rho_J}.$$
Since $\sD_{w_0}\Gro_{ws_i}=\Gro_{w_1}=\sD_{w_0}\Gro_{ws_j}$,
we have
$$\sum_{x\in W_{\{i,j\}}}(-1)^{\ell(x)}x
\e^{\rho_J}\Gro_{ws_i}=
\sum_{x\in W_{\{i,j\}}}(-1)^{\ell(x)}x
\e^{\rho_J}\Gro_{ws_j}.$$
It remains to remark that
\eqn
\sum_{x\in W_{\{i,j\}}}(-1)^{\ell(x)}x
\e^{\rho_J}\Gro_{ws_i}
&=&\sum_{x\in W_{\{i,j\}},\,xs_i>x}(-1)^{\ell(x)}x(1-s_i)
\e^{\rho_J}\Gro_{ws_i}\\
&=&\sum_{x\in W_{\{i,j\}},\,xs_i>x}(-1)^{\ell(x)}x
\e^{\rho_J}(1-\e^{-\al_i})\Gro_{ws_i}.
\eneqn

\medskip\noindent
{\em Proof of \eqref{eq:last}} (iii)\quad
Let $w_0$ be the longest element of $W_{\{i,j\}}$,
and let $w_1$ be the
shortest element of
$wW_{\{i,j\}}$.
Hence $w_1w_0$ is the longest element of $wW_{\{i,j\}}$.

It is enough to show that
$$\sum_{x\in W_{\{i,j\}},\,xs_i>x}(-1)^{\ell(x)}x
\e^{\Lambda_i}(1-\e^{-\alpha_i})\Gro_{ws_i}=0.$$
Since $\e^{\Lambda_i}(1-\e^{-\alpha_i})\Gro_{ws_i}=
(1-s_i)\e^{\Lambda_i}\Gro_{ws_i}$,
it is enough to show
\eq
\sum_{x\in W_{\{i,j\}}}(-1)^{\ell(x)}x\e^{\Lambda_i}\Gro_{ws_i}=0.
\eneq
Since $ws_j>w>ws_i$,
we have $\sD_{w_0s_j}\Gro_{ws_i}=\Gro_{w_1}=\sD_{w_0}\Gro_{ws_i}$.
Hence we have
$$
(\sD_{w_0s_j}-\sD_{w_0})\Gro_{ws_i}=0.
$$
Setting $K=(\sD_{w_0s_j}-\sD_{w_0})\circ \e^{-\Lambda_i}$,
we obtain
\eq
K\e^{\Lambda_i}\Gro_{ws_i}=0.\label{eq:KA}
\eneq
Since
$K=(\sD_{w_0s_j}-\sD_{w_0})\circ\sD_i\circ \e^{-\Lambda_i}$
and $\sD_i\circ \e^{-\Lambda_i}=\e^{-\Lambda_i}(1-\e^{-\al_i})^{-1}(1-s_i)$,
we have
\eq
K\circ(1+s_i)=0.\label{eq:Ki}
\eneq
On the other hand, we have
$K\circ\sD_j=
(\sD_{w_0s_j}-\sD_{w_0})\circ \sD_j\circ \e^{-\Lambda_i}=0$.
Since
$\sD_j=(1+s_j)\circ(1-\e^{-\al_j})^{-1}$, we have
$0=K\circ (1+s_j)\circ (1-\e^{-\al_j})^{-1}$, which implies
$K\circ (1+s_j)=0$.
Together with \eqref{eq:Ki},
$K$ can be written as $K=\psi\circ E$ for some $\psi$
in the quotient field of $\Z[P]$.
Here, $E=\sum_{x\in W_{\{i,j\}}}(-1)^{\ell(x)}x$.
Since $K\not=0$, $\psi$ does not vanish and
\eqref{eq:KA} implies the desired result $E\e^{\Lambda_i}\Gro_{ws_i}=0$.

This completes the proof of \eqref{eq:last}.

\section{Global cohomology character formulas}
\label{sec:charac}

The affine Grothendieck polynomials give the character formula for
the cohomologies of $\O_{X_w}(\mu)\seteq\O_{X_w}\otimes\O_X(\mu)$
under certain conditions on $\mu\in P$.

For $w\in W$, the $B$-orbit
$\dX_w$ is contained in a $T$-stable open affine set
$V_w\seteq wBx_0$ as a closed subset.
As a scheme with $T$-action, $V_w$ is isomorphic to
the group scheme ${}^wU$ whose Lie algebra
is $\soplus_{\al\in w\Delta^+}\g_\al$.
We have a commutative diagram
$$\xymatrix{
{\raisebox{-1ex}[1.5ex][2ex]{$\dX_w$}}\ar@{-}[r]^(.3)\sim\ar@{^{(}->}[d]
&{\raisebox{-.5ex}[1.5ex][2ex]{$\prod_{\al\in\Delta^+\cap w\Delta^+}\g_\al
$}}
\ar@{^{(}->}[d]\\
V_w\ar@{-}[r]^(.35)\sim&{\prod_{\al\in w\Delta^+}\g_\al.}
}$$
Let $\Coh_T(\O_{V_w})$ be the abelian category of
coherent $T$-equivariant $\O_{V_w}$-modules.

\Lemma\label{lem:81}
Any $\F\in \Coh_T(\O_{V_w})$ admits
a free resolution in $\Coh_T(\O_{V_w})$:
\eq&&
0\to F_n\otimes\O_{V_w}\to \cdots\to F_1\otimes\O_{V_w}\to F_0\otimes\O_{V_w}
\to \F\to0,\label{eq:resF}
\eneq
where $F_k$ are finite-dimensional $T$-modules.
\enlemma
\Proof
Set $E=\oplus_{\al\in w\Delta^+}(\g_\al)^*$.
Then $V_w$ is isomorphic to
$\Spec(S(E))$.
Hence, there exists a finite-dimensional
$T$-stable subspace $E'\subset E$ such that
$\F$ is the pull back of a coherent $T$-equivariant
sheaf on $\Spec(S(E'))$ by the faithfully flat projection
$\Spec(S(E))\to\Spec(S(E'))$.
Hence the assertion is a consequence of the following well-known lemma.
\QED

\Lemma
Let $E$ be a finite-dimensional $T$-module
whose weights are contained in $\set{\la\in P}{\lan h,\la\ran>0}$
for some $h\in P^*$.
Then for any $T$-equivariant $\O_E$-module
$\F$, there exists
a free resolution of $\F$ in $\Coh_T(\O_{E})$:
\eqn
0\to F_n\otimes\O_{E}\to \cdots\to F_1\otimes\O_{E}\to F_0\otimes\O_{E}
\to \F\to0,
\eneqn
where $F_k$ are finite-dimensional $T$-modules.
\enlemma

For a locally closed subset $S$ of a topological space $Z$,
we denote by $H^k_S(Z;\sbullet)$ the $k$-th relative cohomology,
and by $\h^k_S(\sbullet)$ the $k$-th local cohomology
(see e.g.\ \cite{H,KS0}).
The following results are proved in \cite{K1}.

\Lemma\label{lem:relc}
For $w\in W$ and $\mu\in P$, we have
\bnum
\item
$H^k_{\dX_w}(X;\O_X(\mu))=0$ for $k\not=\ell(w)$,
\item
$H^{\ell(w)}_{\dX_w}(X;\O_X(\mu))$ is isomorphic to
the dual Verma module with highest weight $w(\mu+\rho)-\rho$, and
$$\ch\bl(H^{\ell(w)}_{\dX_w}(X;\O_X(\mu))\br)
=
\dfrac{\e^{w(\mu+\rho)-\rho}}
{\prod_{\al\in\Delta^+}(1-\e^{-\al})^{\dim \g_\al}}.
$$
\ee
\enlemma
Here, for a $T$-module $M$ such that
its weight space $M_\la$ of weight $\la$ is finite-dimensional
for any $\la\in P$, we set
\eq
&&\ch(M)=\sum_{\la\in P}(\dim M_\la)\,\e^\la.
\eneq

\Lemma
For any $w\in W$ and any coherent $T$-equivariant $\O_{V_w}$-module $\F$,
we have
\bnum
\item
$H^k_{\dX_w}(V_w;\F)=0$ for $k>\ell(w)$, \label{prop1}
\item
For any $\xi\in P$,
$\dim H^k_{\dX_w}(V_w;\F)_\xi<\infty$,\label{prop2}
\item
there exists a finite subset $S$ of $P$ such that
the set of weights of $H^k_{\dX_w}(V_w;\F)$ is contained in $S+Q_-$,
where $Q_-\seteq\sum_i\Z_{\le0}\al_i$,\label{prop3}
\item
$\sum\limits_k(-1)^k\ch(H^k_{\dX_w}(V_w;\F))
=(-1)^{\ell(w)}\ch(H^{\ell(w)}_{\dX_w}(X;\O_X))\cdot
\bl(\sum\limits_k(-1)^k\ch(L_ki_w^*\F)\br)$.\label{prop4} \ee
\enlemma \Proof Since $\dX_w$ is a closed subscheme of the affine
scheme $V_w$ defined as the intersection of the zero loci of $f_i$
($1\le i\le \ell(w)$) for some $f_i\in \O_{V_w}(V_w)$, we obtain
\eqref{prop1}.

Let us prove the other statements.
Let us take a free resolution as in \eqref{eq:resF}.
Then $H^k_{\dX_w}(V_w;\F)$
is the cohomology group of
$F_\bullet\otimes H^{\ell(w)}_{\dX_w}(V_w;\O_X)$.
Hence the results follow from the
corresponding fact for $H^k_{\dX_w}(V_w;\O_X)$  in Lemma~\ref{lem:relc}
and $\sum\limits_k(-1)^k\ch(L_ki_w^*\F)=\sum\limits_k(-1)^k\ch(F_k)$.
\QED

By this lemma, we obtain the following result.

\Prop\label{prop:relc}
\bnum
\item
For $w\in W$ and a $B$-stable quasi-compact open subset $\Omega$ of $X$
such that $\dX_w\subset\Omega$, we have homomorphisms
\eq
&&K_B(X)\To K_{B}(\Omega)\To
\prod_{\la\in P}\Z\,\e^\la
\eneq
given by
$[\F]\mapsto \sum\limits_{k=0}^\infty(-1)^k\ch\bl(H^k_{\dX_w}(\Omega;\F)\br)$.
\item
$\sum\limits_k(-1)^k\ch(H^k_{\dX_w}(X;\F))
=(-1)^{\ell(w)}\ch(H^{\ell(w)}_{\dX_w}(X;\O_X))\cdot
\ch(i_w^*([\F]))$
for any $\F\in \Coh_B(\O_X)$.
\ee
\enprop

\Lemma\label{lem:relv}
Let $w,x\in W$ and $\mu\in P$.
Then, we have
$$\mbox{$\h^k_{\dX_x}(\O_{X_w})=0$
and $H^k_{\dX_x}(X;\O_{X_w}(\mu))=0$
unless $x\ge w$ and $k=\ell(x)-\ell(w)$.}
$$
\enlemma
\Proof
We may assume that  $X_x\subset X_w$.
Set $\omega_{X_w}=\ext^{\ell(w)}_{\O_X}(\O_{X_w},\O_X)$.
Since $\O_{X_w}$ is Cohen-Macaulay by Proposition~\ref{prop:cm},
we have $\O_{X_w}=\Rhom_{\O_X}(\omega_{X_w},\O_X)[\ell(w)]$,
and
$$\h^k_{\dX_x}(\O_{X_w})=\ext^{k+\ell(w)-\ell(x)}_{\O_X}\br(\omega_{X_w},
\h^{\ell(x)}_{\dX_x}(\O_{X})\br).$$
Let $\xi$ be a generic point of $\dX_x$.
Since $\h^{\ell(x)}_{\dX_x}(\O_{X})_\xi$ is
an injective $(\O_X)_\xi$-module (see \cite{H}),
$\h^k_{\dX_x}(\O_{X_w})_\xi=0$ for $k\not=\ell(x)-\ell(w)$.
Since $\h^k_{\dX_x}(\O_{X_w})$ is a quasi-coherent $B$-equivariant
$\O_X$-module
and $\dX_x$ is a $B$-orbit, $\h^k_{\dX_x}(\O_{X_w})\vert_{V_w}=0$.
Let $j\cl V_w\hookrightarrow X$ be the inclusion.
Since $j$ is affine, $\h^k_{\dX_x}(\O_{X_w})=j_*j^{-1}\h^k_{\dX_x}(\O_{X_w})
=0$ for $k\not=\ell(x)-\ell(w)$.

Since $\dX_x$ is affine,
$H^k_{\dX_x}(X;\O_{X_w}(\mu))=\Gamma(\dX_x;\h^k_{\dX_x}(\O_{X_w}(\mu)))=0$
for $k\not=\ell(x)-\ell(w)$.
\QED

Note that Lemmas~\ref{lem:81}--\ref{lem:relv} still hold for any
symmetrizable Kac-Moody Lie algebra $\g$. Now we shall use the fact
that $\g$ is affine.

\Lemma\label{lem:chrel} For $w\in W$, let us write \eq
&&\Gro_w=\sum\limits_{(\la,\al)\in L\times Q}a_{\la,\al}\e^{\la+\al}
\quad\text{with $a_{\la,\al}\in R$.}\label{eq:gro} \eneq Then \eqn
\ch\bl(H^{\ell(x)-\ell(w)}_{\dX_x}(X;\O_{X_w}(\mu))\br) &=&
(-1)^{\ell(w)}\, \dfrac{\e^{x(\mu+\rho)-\rho}j_x(\Gro_w)}
{\prod_{\al\in\Delta^+}(1-\e^{-\al})^{\dim \g_\al}}\\
&=&
(-1)^{\ell(w)}\,
\dfrac{\sum\limits_{(\la,\al)\in L\times Q}a_{\la,\al}\e^{x(\mu+\la+\al+\rho)
-\la-\rho}}
{\prod_{\al\in\Delta^+}(1-\e^{-\al})^{\dim \g_\al}}.\eneqn
\enlemma
\Proof
By Lemma~\ref{lem:relc} and Proposition~\ref{prop:relc},
we have
\eqn
&&(-1)^{\ell(x)-\ell(w)}\ch(H^{\ell(x)-\ell(w)}_{\dX_x}(X;\O_{X_w}(\mu)))\\
&&\hs{15ex}=\sum_k(-1)^k\ch(H^k_{\dX_x}(X;\O_{X_w}(\mu)))\\
&&\hs{15ex}=\sum_k(-1)^k\sum_{(\la,\al)\in L\times Q}a_{\la,\al}\,\e^{-\la}
\ch\bl(H^k_{\dX_x}(X;\O_{X}(\mu+\la+\al))\br)\\
&&\hs{15ex}=(-1)^{\ell(x)}\sum_{(\la,\al)\in L\times Q}a_{\la,\al}\,\e^{-\la}
\ch\bl(H^{\ell(x)}_{\dX_x}(X;\O_{X}(\mu+\la+\al))\br),
\eneqn
which implies the desired result.
\QED

For $\ell\in\Z_{\ge0}$, let $(X_w)_\ell=\cup_{x}X_x$
where $x$ ranges over the elements of $W$ such that $x\ge w$ and
$\ell(x)\ge\ell(w)+\ell$.
Then $\{(X_w)_\ell\}_{\ell\in\Z_{\ge0}}$ is
a decreasing sequence of $B$-stable closed subsets of $X_w$.
Moreover, $(X_w)_\ell\setminus(X_w)_{\ell+1}$ is a disjoint union of
$\dX_x$ where $x$ ranges over the elements of $W$ such that
\eq
&&\mbox{$x\ge w$ and $\ell(x)=\ell(w)+\ell$.}\label{eq:xw}
\eneq
By Lemma~\ref{lem:relv}, we have
$$\mbox{$H^{\;k}_{(X_w)_\ell\setminus(X_w)_{\ell+1}}(X;\O_{X_w}(\mu))
=0$ for $k\not=\ell$,}$$
because it is the direct sum of
$H^k_{\dX_x}(X;\O_{X_w}(\mu))$
where $x\in W$ ranges over the elements satisfying \eqref{eq:xw}.
Hence, by a general argument (see \cite{H,Kumar}),
$H^k(X;\O_{X_w}(\mu))$ is the $k$-th cohomology group
of
\eq
&&\ba{lll}
H^{0}_{X_w\setminus(X_w)_{1}}(X;\O_{X_w}(\mu))
&\to&
H^{1}_{(X_w)_1\setminus(X_w)_{2}}(X;\O_{X_w}(\mu))\\[2ex]
&&\qquad\to
H^{2}_{(X_w)_2\setminus(X_w)_{3}}(X;\O_{X_w}(\mu))
\to\cdots.
\ea
\eneq

\Cor\label{cor:multf}
Let $\mu\in P$ and $w\in W$.
Assume that $\lan c,\mu\ran\ge0$.
Then, for any $\xi\in P$,
$\sum\limits_k\dim
H^{k}(X;\O_{X_w}(\mu))_\xi$
is finite.
\encor
\Proof
Let us set $d(\al)=\sum_im_i$
for $\al=\sum_im_i\al_i\in Q$.
With the notation as in \eqref{eq:gro},
if $a_{\la,\al}\not=0$, then
$\lan c,\la+\al+\rho\ran>0$ by Proposition~\ref{prop:main}.
Hence, we have
$\lan c,\mu+\la+\al+\rho\ran>0$.
Therefore, for any integer $n$, there are only finitely many
$x\in W$ such that
$d\bl(\mu+\la+\al+\rho-x(\mu+\la+\al+\rho)\br)<n$.
Hence,
$\sum\limits_{x\in W}\dim
H^{\ell(x)-\ell(w)}_{\dX_x}(X;\O_{X_w}(\mu))_\xi$
is finite by Lemma~\ref{lem:chrel}.
Since $H^k(X;\O_{X_w}(\mu))$ is a subquotient of
$\hs{-2ex}\soplus_{x\in W,\,\ell(x)=k+\ell(w)}H^{k}_{\dX_x}(X;\O_{X_w}(\mu))$,
we obtain the desired result.
\QED

Thus we obtain the following proposition.
\Prop\label{prop:char}
Let $w\in W$ and $\mu\in P$.
Assume that $\lan c,\mu\ran\ge0$.
Then we have, with the notation \eqref{eq:gro},
\eq
&&
\sum_k(-1)^k\ch\bl(H^k(X;\O_{X_w}(\mu))\br)
=
\sum\limits_{(\la,\al)\in L\times Q}a_{\la,\al}\e^{-\la}\chi_{\mu+\la+\al},
\eneq
where
$\chi_\mu=\dfrac{\sum_{x\in W}(-1)^{\ell(x)}\e^{x(\mu+\rho)-\rho}}
{\prod_{\al\in\Delta^+}(1-\e^{-\al})^{\dim \g_\al}}$.
\enprop
Note that $\sum_k(-1)^k\ch\bl(H^k(X;\O_{X_w}(\mu))\br)$
has a sense by Corollary~\ref{cor:multf}.
\Proof
We have
\eqn
\sum_k(-1)^k\ch\bl(H^k(X;\O_{X_w}(\mu))\br)
&=&\sum_k(-1)^k\ch\bl(H^k_{(X_w)_{k}\setminus(X_w)_{k+1}}(X;\O_{X_w}(\mu))\br)
\\
&=&\sum_{x\in W}(-1)^{\ell(x)-\ell(w)}
\ch\bl(H^{\ell(x)-\ell(w)}_{\dX_x}(X;\O_{X_w}(\mu))\br),
\eneqn
and Lemma~\ref{lem:chrel} implies the desired result.
\QED

\begin{conjecture}\label{conj:van}
We conjecture that,
if $\mu$ is dominant, then $H^k(X;\O_{X_w}(\mu))=0$
for $k\not=0$ and
\eq\Gamma(X;\O_{X}(\mu))\to\Gamma(X;\O_{X_w}(\mu))\label{eq:res}\eneq
is surjective.
\end{conjecture}

Note that $\Gamma(X;\O_{X}(\mu))$ is isomorphic to
the irreducible $\g$-module $V(\mu)$ with highest weight $\mu$, and
the kernel $N$ of \eqref{eq:res} is equal to
$\set{v\in V(\mu)}{U(\gb)v\cap V(\mu)_{w\mu}=0}$.
The module $V(\mu)$ has a non-degenerate symmetric bilinear form
with respect to which the $e_i$'s and the $f_i$'s are adjoint to each other,
and $N$ is orthogonal to $U(\gb_-)u_{w\mu}$.
Here $u_{w\mu}$ is a non-zero vector in the one-dimensional weight space
$V(\mu)_{w\mu}$.
Hence if the conjecture is true,
$\soplus_{\xi\in P}\bl(\Gamma(X;\O_{X_w}(\mu))_\xi\br)^*$
is isomorphic to $U(\gb_-)u_{w\mu}\subset V(\mu)$.
By Proposition~\ref{prop:char}, Conjecture~\ref{conj:van} implies
\eq\quad&&\ch(U(\gb_-)u_{w\mu})=\hs{-2ex}
\sum\limits_{(\la,\al)\in L\times Q}\hs{-1.5ex}
a_{\la,\al}\e^{-\la}\chi_{\mu+\la+\al}
\quad\text{for any dominant integral weight $\mu$.}\eneq

\section{Equivariant cohomology}
\label{sec:eqcoh}

Theorem~\ref{th:main} implies a similar result on the equivariant
cohomology of affine flag manifolds. For a $B$-stable quasi-compact
open subset $\Omega$ of $X$, the equivariant cohomology
$H^*_B(\Omega,\C)$ is a free module over the ring
$H^*_B(\pt,\C)\simeq\C[\gt]=S(\gt^*)$ generated by the
$B$-equivariant cohomology classes $[X_w]$ ($\dX_w\subset\Omega$). For
any $k$, $H_B^k(X,\C)\simeq H^k_B(\Omega,\C)$ if $\Omega$ is large
enough (if $2\,\codim(X\setminus\Omega)>k+1$). Hence we have
$$H^*_B(X,\C)=\soplus_{w\in W}H^*_B(\pt,\C)[X_w].$$
We have a homomorphism
$S(\gt^*)\to H^*_B(X,\C)$
by $P\ni\la\mapsto c_2(\O_{X}(\la))\in H^2_B(X,\C)$,
where $c_2$ is the second Chern class.
It induces an $H^*_B(\pt)$-linear homomorphism
$$H^*_B(\pt)\otimes_{S(\gt^*)^W} S(\gt^*)\to H^*_B(X,\C).$$
Let us write $\widehat{H}^*_B(X,\C)=\prod_kH^k_B(X,\C)$. The equivariant
Chern character defines a homomorphism
$$\ch_B\cl K_B(X)\to \widehat{H}^*_B(X,\C).$$

Hence we have a commutative diagram
$$\xymatrix{
K_B(\pt)\otimes\Z[P]\ar[r]\ar[d]_{\exp\otimes\exp}
&K_B(X)\ar[d]^{{\ch_B}}\\
{\widehat{H}_B^*(\pt,\C)\otimes \widehat{S}(\gt^*)}\ar[r]
&{\widehat{H}_B^*(X,\C)
}}$$
where $\widehat{S}(\gt^*)=\prod_nS^n(\gt^*)$ and
$\exp\cl \Z[P]\to \widehat{S}(\gt^*)$ is given by
$P\ni\la\mapsto \sum_n\la^n/n!\in \widehat{S}(\gt^*)$.

Since the component of $\ch_B(\O_{X_w})$
of degree $2\ell(w)$ coincides with
$[X_w]$, we can translate
Theorem~\ref{th:main} as follows.

\Theorem $\C[\delta,\delta^{-1}]\otimes_{\C[\delta]}
H^*_B(\pt,\C)\otimes_{S(\gt^*)^W} S(\gt^*) \isoto
\C[\delta,\delta^{-1}]\otimes_{\C[\delta]}H^*_B(X,\C).$ \entheorem
Note that $S(\gt^*)^W\simeq S({\gt^*}^W)[\Delta]$, where $\Delta$ is
the Casimir operator
\eqn
\Delta&=&\sum_{i\in I}\dfrac{1}{(\al_i,\al_i)}\La_i\cdot\al_i
-\sum_{i,j\in I}\dfrac{(\al_i,\al_j)}{(\al_i,\al_i)(\al_j,\al_j)}
\La_i\cdot\La_j\\
&=&\dfrac{1}{2}
\sum_{i\in I}\dfrac{1}{(\al_i,\al_i)}\La_i\cdot\bl(\al_i+\eta(\al_i)\br).
\eneqn
\begin{remark}
The absolute $K$-group $K(X)$ is similarly defined. It is isomorphic
to $\prod_{w\in W}\Z[\O_{X_w}]$ and there is a Borel map $\Z[P]\to
K(X)$. However, Theorem~\ref{th:main} gives no information on this
map because $1-\e^{n\delta}$ vanishes in $K(X)$.
\end{remark}

\section{Examples of affine Grothendieck polynomials}
\label{sec:ex}
It is easy to verify directly that
\begin{align*}
  \Gro_{s_i} = 1 - \e^{-\La_i}\qquad\text{for all $i\in I$.}
\end{align*}
Indeed, we have an exact sequence in $\Coh_B(\O_X)$:
$$0\to\C_{\La_i}\otimes\O_X(-\La_i)\to\O_X\to\O_{X_{s_i}}\to 0,$$
where $\C_{\La_i}$ is the one-dimensional $B$-module with weight $\La_i$.

Slightly more generally, if $\set{s_j}{j\in J}$ is a collection of
mutually commuting simple reflections for some $J\subset I$, then
\begin{align*}
  \Gro_{\prod_{j\in J} s_j} = \prod_{j\in J} (1-\e^{-\La_j}).
\end{align*}

The proof of the existence of the elements $\Gro_w$ given in this
paper yields an algorithm to compute them. We have implemented this
algorithm and used it to provide the examples below, using the
notation $q=\e^\delta$ and $\mm^\la=\sum_{\mu\in
W_{I\setminus\{0\}} \cdot \la} \e^\mu$ for the sum of the
exponentials of the weights in the orbit of the element $\la\in P$
under the Weyl group $W_{I\setminus\{0\}}$ of the classical
subalgebra. Usually $\la$ is taken to be anti-dominant with respect
to $W_{I\setminus\{0\}}$. We write only the subscripts of the
simple reflections to indicate Weyl group elements, so that
$\Gro_{10}$ means $\Gro_{s_1s_0}$, for example.

\smallskip\noindent
For $A_{n-1}^{(1)}$ ($n\ge2$) ($I=\Z/n\Z=\{0,1,\ldots,n-1\}$):
\begin{align*}
  \Gro_{10} &= 1+\bl(1-q\br)^{-1}
\bl\{q\mm^{-\La_0}-\mm^{-\La_1}
+\sum_{k=1}^{n-1}\mm^{-\La_1-\La_k+\La_{k+1}+\al_1+\cdots+\al_k} \}.
\end{align*}

\smallskip\noindent
For $A_1^{(1)}$:
{\allowdisplaybreaks
\begin{align*}
\Gro_{010} &= \,1+\bl(\bl(1-q\br)\bl(1-q^2\br)\br)^{-1} \bl\{
-(1+q)(1+q^2)\mm^{-\La_0}
-q\mm^{-\La_0-\al_1} \\*
&\,-q^2\mm^{\La_0-2\La_1+\al_1} +(q+q^2)\mm^{-\La_1}
+(q+q^2)\mm^{-2\La_0+\La_1}-q^2\mm^{-3\La_0+2\La_1-\al_1} \br\},
\\
\Gro_{1010} &= \,1+((1-q)(1-q^2)(1-q^3))^{-1} \bl\{
q(1+q+q^2)(1+q+q^3)\mm^{-\La_0}\\*
&+(q^2+q^3+q^4)\mm^{-\La_0-\al_1}
+q^2\mm^{3\La_0-4\La_1+2\al_1}-(q+q^2+q^3)\mm^{2\La_0-3\La_1+\al_1}
\\*
&+(q+q^2+q^3)\mm^{\La_0-2\La_1}+(1+q+q^2)(1+q+q^3)\mm^{\La_0-2\La_1+\al_1}
\\*
&-(1+q+q^2)(1+q^2+q^3)\mm^{-\La_1}-q^2\mm^{-\La_1-\al_1} \\
&-(q^3+q^4+q^5)\mm^{-2\La_0+\La_1} + q^5\mm^{-3\La_0+2\La_1-\al_1}
 \br\}, \\
 \Gro_{01010} &= 1+((1-q)(1-q^2)(1-q^3)(1-q^4))^{-1}\bl\{ \\*
 &-(1+q^2)(1+q+q^2)(1+q^2+2q^3+q^4+q^6)\mm^{-\La_0}
 -q^4\mm^{-\La_0-2\al_1} \\*
 &-q(1+q)(1+q^2)(1+q^2+q^3+q^4)\mm^{-\La_0-\al_1}
 -q^6\mm^{3\La_0-4\La_1+2\al_1} \\*
 &+q^4(1+q)(1+q^2)\mm^{2\La_0-3\La_1+\al_1}
 -q^3(1+q^2)(1+q+q^2)\mm^{\La_0-2\La_1} \\*
 &-q^2(1+q)(1+q^2)(1+q+q^2+q^4)\mm^{\La_0-2\La_1+\al_1} \\*
 &+q(1+q)^2(1+q^2)(1+q^2+q^4)\mm^{-\La_1}+q^3(1+q)(1+q^2)\mm^{-\La_1-\al_1}
 \\*
 &+q(1+q)^2(1+q^2)(1+q^2+q^4)\mm^{-2\La_0+\La_1}
+q^3(1+q)(1+q^2)\mm^{-2\La_0+\La_1-2\al_1}
 \\*
 &-q^3(1+q^2)(1+q+q^2)\mm^{-3\La_0+2\La_1} \\*
 &-q^2(1+q)(1+q^2)(1+q+q^2+q^4)\mm^{-3\La_0+2\La_1-\al_1}
 \\*
 &+q^4(1+q)(1+q^2)\mm^{-4\La_0+3\La_1-2\al_1} -q^6\mm^{-5\La_0+4\La_1-2\al_1}
 \bl\}.
\end{align*}
}

\smallskip\noindent
For $A_2^{(1)}$:
{\allowdisplaybreaks
\begin{align*}
  \Gro_{010} &=
  1-\e^{-\La_0-\La_1}+\bl(1-q\br)^{-1}\bl\{-q\e^{-\La_0}(\e^{-\al_1}+\e^{-\al_1-\al_2})
  \\*
  &-\e^{-\La_1}(\e^{\al_1}+\e^{\al_1+\al_2})+\e^{\La_0-\La_1-\La_2+\al_1+\al_2}+q\e^{-\La_0+\La_1-\La_2-\al_1} \\
  &+ q\e^{-2\La_0+\La_2-\al_1-\al_2}+\e^{-2\La_1+\La_2+\al_1} \br\}, \\
  \Gro_{210} &=
  1 + ((1-q)(1-q^2))^{-1}
  \bl\{-(q^2+q^3)\mm^{-\La_0}+(q+q^2)\mm^{-\La_1} \\*
  &-\mm^{2\La_0-3\La_2+\al_1+2\al_2}+(1+q)\mm^{\La_0-2\La_2+\al_2}-(q+q^2)\mm^{\La_1-2\La_2+\al_2}
  \\*
  &-(1+q)\mm^{-\La_2}-q\mm^{-\La_2-\al_1}-(q+q^2)\mm^{\La_0-\La_1-\La_2+\al_1+\al_2}
  \\*
  &+(q+q^2)\mm^{-\La_0+\La_1-\La_2-\al_1}-q^2\mm^{-2\La_0+2\La_1-\La_2-\al_1}-q^2\mm^{-2\La_1+\La_2+\al_1}
    \br\}, \\
  \Gro_{1210} &= 1+q^{-1}\mm^{-\La_1-\La_2+\al_1+\al_2} +
  ((1-q)(1-q^2))^{-1}\bl\{ -(1+q+q^2+q^3)\mm^{-\La_0} \\*
  &-\mm^{2\La_0-3\La_1+2\al_1+\al_2}+(1+q)\mm^{\La_0-2\La_1+\al_1}-q\mm^{-\La_1-\al_2-\al_3}
  -\mm^{2\La_0-3\La_2+\al_1+2\al_2}
  \\*
  &+(1+q)\mm^{\La_0-2\La_2+\al_2}-(1+q+q^2)\mm^{\La_1-2\La_2+\al_2}-q\mm^{-\La_2-\al_1}
  \\*
  &-q^{-1}(1+q+q^2+q^3)\mm^{\La_0-\La_1-\La_2+\al_1+\al_2}+(q+q^2)\mm^{-\La_0+\La_1-\La_2-\al_1}
  \\*
  &
  -q^2\mm^{-2\La_0+2\La_1-\La_2-\al_1}-(1+q+q^2)\mm^{-2\La_1+\La_2+\al_1}\\*
  &+(q+q^2)\mm^{-\La_0-\La_1+\La_2-\al_2-\al_3}  -q^2\mm^{-2\La_0-\La_1+2\La_2-\al_1-2\al_2-\al_3}
  \br\}.
\end{align*}
}
\smallskip\noindent
For $A_3^{(1)}$:
\begin{align*}
  \Gro_{210} &= 1 +((1-q)(1-q^2))^{-1}
  \bl\{-(q^2+q^3)\mm^{-\La_0}+(q+q^2)\mm^{-\La_1} \\
  &-(1+q)\mm^{\La_0-2\La_2+\al_1+2\al_2+\al_3}
  -(1+q)\mm^{-\La_2}-q\mm^{-\La_2-\al_1} \\
  &+(q+q^2)\mm^{-\La_0+\La_1-\La_2-\al_1}-q^2\mm^{-2\La_0+2\La_1-\La_2-\al_1}-q^2\mm^{-2\La_1-\La_2+\al_1}\\
  &-\mm^{2\La_0-\La_2-2\La_3+\al_1+2\al_2+2\al_3}-(q+q^2)\mm^{\La_0-\La_1-\La_3+\al_1+\al_2+\al_3}
  \\
  &+(1+q)\mm^{\La_0-\La_2-\La_3+\al_2+\al_3}-(q+q^2)\mm^{\La_1-\La_2-\La_3+\al_2+\al_3}
  \\
  &+(1+q)\mm^{-2\La_2+\La_3+\al_2}-(q+q^2)\mm^{-\La_0+\La_1-2\La_2+\La_3+\al_2}
  \\
  &-(q+q^2)\mm^{-\La_1-\La_2+\La_3+\al_1+\al_2}-\mm^{-3\La_2+2\La_3+\al_1+2\al_2}
  \br\}.
\end{align*}

\smallskip\noindent
For $D_4^{(1)}$ ($I=\{0,1,2,3,4\}$ where $2$ is the node connected
to all others):
{\allowdisplaybreaks
\begin{align*}
  \Gro_{12} &= (1-\e^{-\La_1})(1-\e^{-\La_1+\al_1}), \\
  \Gro_{121} &=
  (1-\e^{-\La_1})(1-\e^{-\La_1+\al_1})(1-\e^{\La_1-\La_2-\al_1}), \\
  \Gro_{321} &= 1 - \e^{-2\La_3-\al_1+\al_3} + (1-q)^{-1} \bl\{
q\e^{-\La_1}+q\e^{\La_1-\La_2-\al_1}+\e^{\La_1-2\La_3-\al_1+\al_3}
\\*
&+\,\e^{-\La_1+\La_2-2\La_3+\al_3} -
\e^{-\La_3}-q\e^{-\La_3-\al_1}-q\e^{-\La_3-\al_1-\al_2}
-\e^{-\La_3+\al_3}\\*
&-\,q\e^{-\La_3-\al_1+\al_3} - \e^{-\La_3+\al_2+\al_3} - q
\e^{-\La_3-\al_1-\al_2-\al_4}-\e^{-\La_3+\al_2+\al_3+\al_4}\\*
&+\,\e^{\La_0-\La_3-\La_4+\al_2+\al_3+\al_4} + q
\e^{-\La_0+\La_2-\La_3-\La_4-\al_1-\al_2} \\*
&+\, q\e^{-\La_0-\La_3+\La_4-\al_1-\al_2-\al_4} +
\e^{\La_0-\La_2-\La_3+\La_4+\al_2+\al_3}
  \bl\}.
\end{align*}
}

\smallskip\noindent
For $C_n^{(1)}$ let $I=\{0,1,\dotsc,n\}$ and let $\epsilon_i$ be the
$i$-th standard basis vector in $\Z^n$ for $1\le i\le n$ and let
$\al_i=\epsilon_i-\epsilon_{i+1}$ for $1\le i\le n-1$ and
$\al_n=2\epsilon_n$. Then we have
\begin{align*}
  \Gro_{10} &= 1+(1-q)^{-1}\bl\{ q\mm^{-\La_0}-\mm^{-\La_1} \\
  &
  +\sum_{k=1}^n \mm^{-\La_1+\La_{k-1}-\La_k+\epsilon_1+\epsilon_k} +
  \sum_{k=1}^{n-1}\mm^{-\La_1-\La_k+\La_{k+1}+\epsilon_1-\epsilon_{k+1}}\br\}.
\end{align*}

\smallskip\noindent
For $C_2^{(1)}$:
\begin{align*}
\Gro_{010} &=1+((1-q)(1-q^2))^{-1}\bl\{ -(1+q)^3 \mm^{-\La_0}-(q+q^2)\mm^{-\La_0-\al_1}-q\mm^{-\La_0-\al_2} \\
&-q^2\mm^{\La_0-2\La_1+2\al_1+\al_2}+(q+q^2)\mm^{-\La_1}+(q+q^2)\mm^{-2\La_0+\La_1-2\al_1-\al_2} \\
&-q^2\mm^{-3\La_0+2\La_1-2\al_1-\al_2}-q^2\mm^{-\La_0+2\La_1-2\La_2+\al_2}-(q+q^2)\mm^{-\La_2+\al_1+\al_2} \\
&+(q+q^2)\mm^{-\La_0+\La_1-\La_2-\al_1}-\!(q+q^2)\mm^{-2\La_0+2\La_1-\La_2-\al_1}-\!(q+q^2)\mm^{-2\La_0+\La_2-\al_1-\al_2}
\\ &-(q+q^2)\mm^{-2\La_1+\La_2+\al_1}+(q+q^2)\mm^{-\La_0-\La_1+\La_2-\al_1-\al_2}-q^2\mm^{-\La_0-2\La_1+2\La_2-\al_2}
\br\}.
\end{align*}

\section{Proof of Lemma~\ref{lem:inj}}
\label{sec:inj}

In this section, we shall give a proof of  Lemma~\ref{lem:inj}.
Let us set $\gt^*_\cla=\gt^*/(\gt^*)^W$ and
let $\cla\cl \gt^*\to\gt^*_\cla$ be the canonical projection.
Let us recall that
$Q_\cla=Q/\Z\delta=\cla(Q)$ and $W_\cla$ is the image of $W$ in $\Aut(Q_\cla)$.
Note that $\C\otimes Q_\cla$ is a hyperplane of $\gt^*_\cla$.

For $\xi_0\in \C\otimes Q_\cla$, we define $t(\xi_0)\in\Aut(\gt^*)$
by
\eq
&&t(\xi_0)(\la)=\la+(\la,\delta)\xi-(\la,\xi)\delta
-\dfrac{(\la,\delta)(\xi,\xi)}{2}\delta
\eneq
for $\xi\in\C\otimes Q$ such that $\cla(\xi)=\xi_0$.
It does not depend on the choice of $\xi$.

Let $\tQ$ be the sublattice
$Q_\cla\cap\cla\bl(\soplus_i\frac{2}{(\al_i,\al_i)}\al_i\br)$.
Then we have an exact sequence
$$1\To\tQ\To[\ t\ ] W\To W_\cla\To1.$$

\smallskip
Now we shall prove Lemma~\ref{lem:inj}.
Assume that $u\in\Z[P]\otimes_{\Z[P^W]}\Z[P]\simeq\Z[L]\otimes_\Z\Z[P]$
satisfies $i_w^*\beta(u)=0$
for all $w\in W$.
We can write uniquely
$$\mbox{$u=\sum_{\la\in L}(\e^\la\otimes 1)\cdot\xi(u_\la)$
with $u_\la\in\Z[P]$.}$$
Then we have
$\sum_{\la\in L}\e^\la j_w(u_\la)=0$ for any $w\in W$
(see \eqref{def:jw} and the commutative diagram \eqref{com:jw}).
Since $j_w(u_\la)\in\Z[Q]$, we have
$j_w(u_\la)=0$ for any $\la\in L$ and $w\in W$.
Hence we have reduced Lemma~\ref{lem:inj} to the following lemma.
\Lemma
Let $u\in\Z[P]$.
If
$j_{t(\xi_0)}(u)=0$ for any $\xi_0\in\tQ$,
then $u=0$.
\enlemma
\Proof
Write $u=\sum_{(\la,\al)\in L\times Q}a_{\la,\al}\e^{\la+\al}$.
For any $\xi\in Q\cap\cla^{-1}(\tQ\setminus\{0\})$ and any integer $n$,
we have $(\xi,\xi)>0$, and $j_{t(\cla(n\xi))}(u)=0$ reads as
$$\sum_{(\la,\al)\in L\times Q}a_{\la,\al}\e^%
{\al+n\{(\la,\delta)\xi-(\la+\al,\xi)\delta\}
-n^2\frac{(\la,\delta)(\xi,\xi)}{2}\delta}
=0.$$
Hence we have for all $n$ and $\ell$
$$\sum_{(\la,\al)\in L\times Q,\,(\la,\de)=\ell}
a_{\la,\al}\e^%
{\al+n\{(\la,\delta)\xi-(\la+\al,\xi)\delta\}}
=0,$$
which implies that
$$\sum_{(\la,\al)\in L\times Q,\,(\la,\de)=\ell}
a_{\la,\al}\e^%
{\al-n(\la+\al,\xi)\delta}
=0.$$
Hence for all $\ell$, $m$ and $\xi$, we have
$$\sum_{(\la,\al)\in L\times Q,\, (\la,\de)=\ell,\,(\la+\al,\xi)=m}
a_{\la,\al}\e^\al
=0,$$
which implies that
$$\sum_{\la\in L,\,(\la,\xi)=m}
a_{\la,\al}
=0$$
for any $\al$.

Set $\supp(u)=\set{\la\in L}{\text{$a_{\la,\al}\not=0$ for some $\al\in Q$}}$.
Since there exists $\xi\in Q$ such that
$(\la,\xi)\not=(\la',\xi)$ for any pair of distinct
elements $\la$, $\la'$ in $\supp(u)$,
we have
$a_{\la,\al}=0$ for all $\la$, $\al$.
\QED

\begin{remark}\label{rem:inf}
If $\la\in P\setminus P^W$, then $W\la$ is an infinite set.
Indeed, there exists
$\xi\in Q\cap\cla^{-1}(\tQ\setminus\{0\})$ such that $(\la,\xi)\not=0$,
and we can easily see that
$\set{t(\cla(n\xi))\la}{n\in\Z}$ is an infinite set.
\end{remark}

\end{document}